\documentclass[oneside,12pt]{article}
\usepackage[colorlinks,
            citecolor=DarkRed,
            urlcolor=DarkRed]{hyperref} 
\usepackage{tocloft}                    
\usepackage{indentfirst}
\usepackage{amsmath}
\usepackage{amsfonts}
\usepackage{amssymb}
\usepackage{pict2e}
\usepackage{stmaryrd}
\usepackage{bm}                        
\usepackage[x11names,svgnames]{xcolor} 
\usepackage{proof}   
\input diagxy        
\xyoption{curve}     
\usepackage{ifluatex}             %
\ifluatex                         %
  \input edrxaccents.tex          
  \input edrx21chars.tex          
  \input edrxheadfoot.tex         
\else                             %
  \usepackage[utf8]{inputenc}     %
\DeclareUnicodeCharacter{00A1}{\text{\textexclamdown}} 
\DeclareUnicodeCharacter{00A7}{\S}                 
\DeclareUnicodeCharacter{00AC}{\neg}               
\DeclareUnicodeCharacter{00B0}{^\circ}             
\DeclareUnicodeCharacter{00B1}{\pm}                
\DeclareUnicodeCharacter{00B2}{\superscriptTwo}    
\DeclareUnicodeCharacter{00B3}{\superscriptThree}  
\DeclareUnicodeCharacter{00B7}{\cdot}              
\DeclareUnicodeCharacter{00B9}{\superscriptOne}    
\DeclareUnicodeCharacter{00D7}{\times}             
\DeclareUnicodeCharacter{00F7}{\div}               
\DeclareUnicodeCharacter{0393}{\Gamma}             
\DeclareUnicodeCharacter{0394}{\Delta}             
\DeclareUnicodeCharacter{0398}{\Theta}             
\DeclareUnicodeCharacter{039B}{\Lambda}            
\DeclareUnicodeCharacter{03A0}{\Pi}                
\DeclareUnicodeCharacter{03A3}{\Sigma}             
\DeclareUnicodeCharacter{03A6}{\Phi}               
\DeclareUnicodeCharacter{03A8}{\Psi}               
\DeclareUnicodeCharacter{03A9}{\Omega}             
\DeclareUnicodeCharacter{03B1}{\alpha}             
\DeclareUnicodeCharacter{03B2}{\beta}              
\DeclareUnicodeCharacter{03B3}{\gamma}             
\DeclareUnicodeCharacter{03B4}{\delta}             
\DeclareUnicodeCharacter{03B5}{\varepsilon}        
\DeclareUnicodeCharacter{03B6}{\zeta}              
\DeclareUnicodeCharacter{03B7}{\eta}               
\DeclareUnicodeCharacter{03B8}{\theta}             
\DeclareUnicodeCharacter{03B9}{\iota}              
\DeclareUnicodeCharacter{03BA}{\kappa}             
\DeclareUnicodeCharacter{03BB}{\lambda}            
\DeclareUnicodeCharacter{03BC}{\mu}                
\DeclareUnicodeCharacter{03BD}{\nu}                
\DeclareUnicodeCharacter{03BE}{\xi}                
\DeclareUnicodeCharacter{03C0}{\pi}                
\DeclareUnicodeCharacter{03C1}{\rho}               
\DeclareUnicodeCharacter{03C3}{\sigma}             
\DeclareUnicodeCharacter{03C4}{\tau}               
\DeclareUnicodeCharacter{03C6}{\varphi}            
\DeclareUnicodeCharacter{03C7}{\chi}               
\DeclareUnicodeCharacter{03C8}{\psi}               
\DeclareUnicodeCharacter{03C9}{\omega}             
\DeclareUnicodeCharacter{03D5}{\phi}               
\DeclareUnicodeCharacter{1D400}{\mathbf{A}}        
\DeclareUnicodeCharacter{1D401}{\mathbf{B}}        
\DeclareUnicodeCharacter{1D402}{\mathbf{C}}        
\DeclareUnicodeCharacter{1D403}{\mathbf{D}}        
\DeclareUnicodeCharacter{1D404}{\mathbf{E}}        
\DeclareUnicodeCharacter{1D405}{\mathbf{F}}        
\DeclareUnicodeCharacter{1D406}{\mathbf{G}}        
\DeclareUnicodeCharacter{1D407}{\mathbf{H}}        
\DeclareUnicodeCharacter{1D408}{\mathbf{I}}        
\DeclareUnicodeCharacter{1D409}{\mathbf{J}}        
\DeclareUnicodeCharacter{1D40A}{\mathbf{K}}        
\DeclareUnicodeCharacter{1D40B}{\mathbf{L}}        
\DeclareUnicodeCharacter{1D40C}{\mathbf{M}}        
\DeclareUnicodeCharacter{1D40D}{\mathbf{N}}        
\DeclareUnicodeCharacter{1D40E}{\mathbf{O}}        
\DeclareUnicodeCharacter{1D40F}{\mathbf{P}}        
\DeclareUnicodeCharacter{1D410}{\mathbf{Q}}        
\DeclareUnicodeCharacter{1D411}{\mathbf{R}}        
\DeclareUnicodeCharacter{1D412}{\mathbf{S}}        
\DeclareUnicodeCharacter{1D413}{\mathbf{T}}        
\DeclareUnicodeCharacter{1D414}{\mathbf{U}}        
\DeclareUnicodeCharacter{1D415}{\mathbf{V}}        
\DeclareUnicodeCharacter{1D416}{\mathbf{W}}        
\DeclareUnicodeCharacter{1D417}{\mathbf{X}}        
\DeclareUnicodeCharacter{1D418}{\mathbf{Y}}        
\DeclareUnicodeCharacter{1D419}{\mathbf{Z}}        
\DeclareUnicodeCharacter{1D41B}{\mathbf}           
\DeclareUnicodeCharacter{1D422}{\textsl}           
\DeclareUnicodeCharacter{1D42B}{\mathrm}           
\DeclareUnicodeCharacter{1D42C}{\mathsf}           
\DeclareUnicodeCharacter{1D42D}{\text}             
\DeclareUnicodeCharacter{1D4D0}{\mathcal{A}}       
\DeclareUnicodeCharacter{1D4D1}{\mathcal{B}}       
\DeclareUnicodeCharacter{1D4D2}{\mathcal{C}}       
\DeclareUnicodeCharacter{1D4D3}{\mathcal{D}}       
\DeclareUnicodeCharacter{1D4D4}{\mathcal{E}}       
\DeclareUnicodeCharacter{1D4D5}{\mathcal{F}}       
\DeclareUnicodeCharacter{1D4D6}{\mathcal{G}}       
\DeclareUnicodeCharacter{1D4D7}{\mathcal{H}}       
\DeclareUnicodeCharacter{1D4D8}{\mathcal{I}}       
\DeclareUnicodeCharacter{1D4D9}{\mathcal{J}}       
\DeclareUnicodeCharacter{1D4DA}{\mathcal{K}}       
\DeclareUnicodeCharacter{1D4DB}{\mathcal{L}}       
\DeclareUnicodeCharacter{1D4DC}{\mathcal{M}}       
\DeclareUnicodeCharacter{1D4DD}{\mathcal{N}}       
\DeclareUnicodeCharacter{1D4DE}{\mathcal{O}}       
\DeclareUnicodeCharacter{1D4DF}{\mathcal{P}}       
\DeclareUnicodeCharacter{1D4E0}{\mathcal{Q}}       
\DeclareUnicodeCharacter{1D4E1}{\mathcal{R}}       
\DeclareUnicodeCharacter{1D4E2}{\mathcal{S}}       
\DeclareUnicodeCharacter{1D4E3}{\mathcal{T}}       
\DeclareUnicodeCharacter{1D4E4}{\mathcal{U}}       
\DeclareUnicodeCharacter{1D4E5}{\mathcal{V}}       
\DeclareUnicodeCharacter{1D4E6}{\mathcal{W}}       
\DeclareUnicodeCharacter{1D4E7}{\mathcal{X}}       
\DeclareUnicodeCharacter{1D4E8}{\mathcal{Y}}       
\DeclareUnicodeCharacter{1D4E9}{\mathcal{Z}}       
\DeclareUnicodeCharacter{2022}{\bullet}            
\DeclareUnicodeCharacter{2026}{\ldots}             
\DeclareUnicodeCharacter{2080}{{{}_0}}             
\DeclareUnicodeCharacter{2081}{{{}_1}}             
\DeclareUnicodeCharacter{2082}{{{}_2}}             
\DeclareUnicodeCharacter{2083}{{{}_3}}             
\DeclareUnicodeCharacter{2084}{{{}_4}}             
\DeclareUnicodeCharacter{2085}{{{}_5}}             
\DeclareUnicodeCharacter{2086}{{{}_6}}             
\DeclareUnicodeCharacter{2087}{{{}_7}}             
\DeclareUnicodeCharacter{2088}{{{}_8}}             
\DeclareUnicodeCharacter{2089}{{{}_9}}             
\DeclareUnicodeCharacter{208A}{{{}_+}}             
\DeclareUnicodeCharacter{208B}{{{}_-}}             
\DeclareUnicodeCharacter{208C}{{{}_=}}             
\DeclareUnicodeCharacter{208D}{{{}_(}}             
\DeclareUnicodeCharacter{208E}{{{}_)}}             
\DeclareUnicodeCharacter{2090}{{{}_a}}             
\DeclareUnicodeCharacter{2091}{{{}_e}}             
\DeclareUnicodeCharacter{2092}{{{}_o}}             
\DeclareUnicodeCharacter{2093}{{{}_x}}             
\DeclareUnicodeCharacter{2095}{{{}_h}}             
\DeclareUnicodeCharacter{2096}{{{}_k}}             
\DeclareUnicodeCharacter{2097}{{{}_l}}             
\DeclareUnicodeCharacter{2098}{{{}_m}}             
\DeclareUnicodeCharacter{2099}{{{}_n}}             
\DeclareUnicodeCharacter{209A}{{{}_p}}             
\DeclareUnicodeCharacter{209B}{{{}_s}}             
\DeclareUnicodeCharacter{209C}{{{}_t}}             
\DeclareUnicodeCharacter{2113}{\ell}               
\DeclareUnicodeCharacter{214B}{\bindnasrepma}      
\DeclareUnicodeCharacter{2190}{\ot}                
\DeclareUnicodeCharacter{2191}{\upto}              
\DeclareUnicodeCharacter{2192}{\to}                
\DeclareUnicodeCharacter{2193}{\dnto}              
\DeclareUnicodeCharacter{2194}{\bij}               
\DeclareUnicodeCharacter{2195}{\updownarrow}       
\DeclareUnicodeCharacter{2196}{\nwarrow}           
\DeclareUnicodeCharacter{2197}{\nearrow}           
\DeclareUnicodeCharacter{2198}{\searrow}           
\DeclareUnicodeCharacter{2199}{\swarrow}           
\DeclareUnicodeCharacter{21A3}{\epito}             
\DeclareUnicodeCharacter{21A4}{\mapsfrom}          
\DeclareUnicodeCharacter{21A6}{\mapsto}            
\DeclareUnicodeCharacter{21C0}{\rightharpoonup}    
\DeclareUnicodeCharacter{21D0}{\Leftarrow}         
\DeclareUnicodeCharacter{21D2}{\funto}             
\DeclareUnicodeCharacter{21D4}{\Leftrightarrow}    
\DeclareUnicodeCharacter{21DD}{\rightsquigarrow}   
\DeclareUnicodeCharacter{2200}{\forall}            
\DeclareUnicodeCharacter{2202}{\partial}           
\DeclareUnicodeCharacter{2203}{\exists}            
\DeclareUnicodeCharacter{2205}{\emptyset}          
\DeclareUnicodeCharacter{2207}{\nabla}             
\DeclareUnicodeCharacter{2208}{\in}                
\DeclareUnicodeCharacter{2216}{\backslash}         
\DeclareUnicodeCharacter{2218}{\circ}              
\DeclareUnicodeCharacter{221A}{\sqrt}              
\DeclareUnicodeCharacter{221E}{\infty}             
\DeclareUnicodeCharacter{2227}{\land}              
\DeclareUnicodeCharacter{2228}{\lor}               
\DeclareUnicodeCharacter{2229}{\cap}               
\DeclareUnicodeCharacter{222A}{\cup}               
\DeclareUnicodeCharacter{222B}{\int}               
\DeclareUnicodeCharacter{223C}{\sim}               
\DeclareUnicodeCharacter{2243}{\simeq}             
\DeclareUnicodeCharacter{2245}{\cong}              
\DeclareUnicodeCharacter{2248}{\approx}            
\DeclareUnicodeCharacter{2260}{\neq}               
\DeclareUnicodeCharacter{2261}{\equiv}             
\DeclareUnicodeCharacter{2264}{\le}                
\DeclareUnicodeCharacter{2265}{\ge}                
\DeclareUnicodeCharacter{2282}{\subset}            
\DeclareUnicodeCharacter{2283}{\supset}            
\DeclareUnicodeCharacter{2286}{\subseteq}          
\DeclareUnicodeCharacter{2287}{\supseteq}          
\DeclareUnicodeCharacter{2293}{\sqcap}             
\DeclareUnicodeCharacter{2294}{\sqcup}             
\DeclareUnicodeCharacter{2295}{\oplus}             
\DeclareUnicodeCharacter{2296}{\ominus}            
\DeclareUnicodeCharacter{2297}{\otimes}            
\DeclareUnicodeCharacter{2298}{\oslash}            
\DeclareUnicodeCharacter{2299}{\odot}              
\DeclareUnicodeCharacter{22A2}{\vdash}             
\DeclareUnicodeCharacter{22A3}{\dashv}             
\DeclareUnicodeCharacter{22A4}{\top}               
\DeclareUnicodeCharacter{22A5}{\bot}               
\DeclareUnicodeCharacter{22A8}{\vDash}             
\DeclareUnicodeCharacter{22C0}{\bigwedge}          
\DeclareUnicodeCharacter{22C1}{\bigvee}            
\DeclareUnicodeCharacter{22C2}{\bigcap}            
\DeclareUnicodeCharacter{22C3}{\bigcup}            
\DeclareUnicodeCharacter{22C4}{\lozenge}           
\DeclareUnicodeCharacter{22C5}{\Box}               
\DeclareUnicodeCharacter{22EE}{\vdots}             
\DeclareUnicodeCharacter{22EF}{\cdots}             
\DeclareUnicodeCharacter{2329}{\langle}            
\DeclareUnicodeCharacter{232A}{\rangle}            
\DeclareUnicodeCharacter{2581}{\_}                 
\DeclareUnicodeCharacter{25A1}{\Box}               
\DeclareUnicodeCharacter{25FB}{\Box}               
\DeclareUnicodeCharacter{266D}{\flat}              
\DeclareUnicodeCharacter{266E}{\natural}           
\DeclareUnicodeCharacter{266F}{\sharp}             
\DeclareUnicodeCharacter{27E6}{\llbracket}         
\DeclareUnicodeCharacter{27E7}{\rrbracket}         
\DeclareUnicodeCharacter{2806}{{:}}                

\fi                               %
\usepackage{edrx21}               
\usepackage[backend=biber,
   style=alphabetic]{biblatex}    
\begin{document}

\ifluatex
  \catcode`\^^J=10
  \directlua{adddednatlualoader = function () DednatRequire.usenewrequire() end}
  \directlua{dofile "dednat6load.lua"}
\else
  \input\jobname.dnt   
  \def\pu{}
\fi

\def\co#1{{%
  \def\%{\char37}%
  \def\\{\char92}%
  \def\^{\char94}%
  \def\~{\char126}%
  \tt#1%
  }}
\def\qco#1{`\co{#1}'}
\def\qqco#1{``\co{#1}''}

\def\ph{\phantom}

\def\respcomp{\mathsf{respcomp}}
\def\respids {\mathsf{respids}}
\def\sqcond  {\mathsf{sqcond}}
\def\assoc   {\mathsf{assoc}}
\def\idL     {\mathsf{idL}}
\def\idR     {\mathsf{idR}}
\def\univ    {\mathsf{univ}}
\def\Ran     {\mathsf{Ran}}

\def\sfC  {\mathsf{C}}
\def\sfSet{\mathsf{Set}}
\def\Ring {\mathbf{Ring}}
\def\nameof#1{\ulcorner#1\urcorner}
\def\catK {\mathbf{K}}
\def\Dn   {\Downarrow}

\def\veq{\rotatebox{90}{$=$}}
\def\liml{\underleftarrow {\lim}{}}
\def\limr{\underrightarrow{\lim}{}}
\def\veq{\rotatebox{90}{$=$}}

\def\Yzero    {\mathsf{Y0}}
\def\Yzeroplus{\mathsf{Y0^+}}
\def\Yone     {\mathsf{Y1}}
\def\Ytwo     {\mathsf{Y2}}
\def\Ythree   {\mathsf{Y3}}
\def\Yfour    {\mathsf{Y4}}
\def\Yfive    {\mathsf{Y5}}

\def\origphi{\phi}

\def\AProofOf   #1{\llangle#1\rrangle}
\def\AllProofsOf#1{\llbracket#1\rrbracket}

\def\DONE{(DONE)}
\def\DONE{}

\def\slowdiag#1{\standout{#1}}
\def\slowdiag#1{\diag{#1}}


%

\title{On the the missing diagrams \\ in Category Theory \\ (first-person version)}

\author{%
  Eduardo Ochs%
  \thanks{eduardoochs@gmail.com}\\
  }

\maketitle


\begin{abstract}

  Most texts on Category Theory are written in a very terse style, in
  which people pretend a) that all concepts are visualizable, and b)
  that the readers can reconstruct the diagrams that the authors had
  in mind based on only the most essential cues. As an outsider I
  spent years believing that the techniques for drawing diagrams were
  part of the oral culture of the field, and that the insiders could
  read texts on CT reconstructing the ``missing diagrams'' in them
  line by line and paragraph by paragraph, and drawing for each page
  of text a page of diagrams with all the diagrams that the authors
  had omitted. My belief was wrong: there are lots of conventions for
  drawing diagrams scattered through the literature, but that unified
  diagrammatic language did not exist. In this chapter I will show an
  attempt to reconstruct that (imaginary) language for missing
  diagrams: we will see an extensible diagrammatic language, called
  DL, that follows the conventions of the diagrams in the literature
  of CT whenever possible and that seems to be adequate for drawing
  ``missing diagrams'' for Category Theory. Our examples include the
  ``missing diagrams'' for adjunctions, for the Yoneda Lemma, for Kan
  extensions, and for geometric morphisms, and how to formalize them
  in Agda.

\end{abstract}

\newpage

%

\renewcommand{\cfttoctitlefont}{\bfseries}
\setlength{\cftbeforesecskip}{2.5pt}

\tableofcontents


%

\section{Introduction}
\label{missing-diagrams}
\label{introduction}

One of the main themes of this text is a diagrammatic language ---
let's call it DL --- that can be used to draw ``missing diagrams'' for
Category Theory. DL is a {\sl reconstructed language}, and it's easier
to explain it if I explain how it was reconstructed, and which of its
conventions were improvised. It is easier to do it in the first
person.

Suppose that your native language is $A$ and you are learning a
language $B$ by a method that includes conversation classes. You will
have to improvise a lot, but you will usually get feedback quickly.
Now suppose that you are studying a language $C$ --- for example,
Aeolic Greek (\cite{CarsonSappho}) --- mostly by yourself, and the
corpus of texts in $C$ is small. A good exercise is to try to write
your thoughts in $C$, using loanwords and improvised syntactical
constructs when needed, but marking mentally the places in which you
had to improvise. In most cases, but not all, you will eventually find
ways to rewrite those parts to make them look more like $C$.

The conventions of DL are explained in sec.\ref{the-conventions}. A
few of them don't correspond to anything that I could find in the
literature; they are listed at the end of that section.

\msk

{\bf This text has two versions.} The title of this one is ``On the
missing diagrams in Category Theory (first-person version)'', and the
title of the other is just ``On the missing diagrams in Category
Theory''; they (will) have the same technical content, but different
styles. Here's how this happened. Two of the editors of the ``Handbook
of Abductive Cognition'' --- Gianluca Caterina and Rocco Gangle ---
asked me if I could rewrite my notes in \cite{FavC} and submit the
rewritten version as a chapter for the Handbook. I did that, and I
wrote this version, that uses the first person a lot --- for example
to tell stories, and to explain why the current version of DL has a
certain convention, that replaces another, older, convention, that
didn't work very well. Then the other editors asked me of I could
prepare another version that would follow the editorial guidelines
more closely, including the part that says to avoid the first person.
At this moment (april/2022) the ``non-first-person version'' is not
ready yet; the ``first-person version'' will probably be more fun to
read.

\bsk

The best way to introduce DL is to tell the full story of how it
evolved from a long sequence of wrong assumptions and from some
``thoughts that I wanted to express in DL''.



Let me start with some quotes. This one is from Eilenberg and Steenrod
(\cite[p.ix]{EilenbergSteenrod}, but I learned it from
\cite[pp.82--83]{Kromer}):

\begin{quotation}

  The diagrams incorporate a large amount of information. Their use
  provides extensive savings in space and in mental effort. In the
  case of many theorems, the setting up of the correct diagram is the
  major part of the proof. We therefore urge that the reader stop at
  the end of each theorem and attempt to construct for himself the
  relevant diagram before examining the one which is given in the
  text. Once this is done, the subsequent demonstration can be
  followed more readily; in fact, the reader can usually supply it
  himself.

\end{quotation}

I spent a {\sl lot} of my time studying Category Theory trying to
``supply the diagrams myself''. In \cite{EilenbergSteenrod} supplying
the diagrams is not very hard (I guess), but in books like
\cite{CWM2}, in which most important concepts involve several
categories, I had to rearrange my diagrams hundreds of times until I
reached ``good'' diagrams...



The problem is that I expected too much from ``good'' diagrams. The
next quotes are from the sections 1 and 12 of an article that I wrote
about that (\cite{IDARCT}):

\begin{quotation}

  My memory is limited, and not very dependable: I often have to
  rededuce results to be sure of them, and I have to make them fit in
  as little ``mental space'' as possible...

  Different people have different measures for ``mental space'';
  someone with a good algebraic memory may feel that an expression
  like $\mathsf{Frob}: Σ_f(P ∧ f^* Q) ≅ Σ_f P ∧ Q$ is easy to
  remember, while I always think diagramatically, and so what I do is
  that I remember this diagram,

  \begin{center}
  \includegraphics[height=60pt]{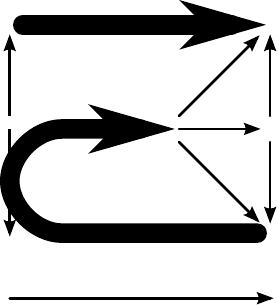}
  \end{center}

  \noindent and I reconstruct the formula from it.

\end{quotation}

\begin{quotation}

  Let's call the ``projected'' version of a mathematical object its
  ``skeleton''. The underlying idea in this paper is that for the
  right kinds of projections, and for some kinds of mathetical
  objects, it should be possible to reconstruct enough of the original
  object from its skeleton and few extra clues --- just like
  paleontologists can reconstruct from a fossil skeleton the look of
  an animal when it was alive.

\end{quotation}

I was searching for a diagrammatic language that would let me express
the ``skeletons'' of categorical definitions and proofs. I wanted
these skeletons to be easy to remember --- partly because they would
have shapes that were easy to remember, and partly because they would
be similar to ``archetypal cases'' (\cite[section 16]{IDARCT}).

\bsk


In 2016 and 2017 I taught a seminar course for undergraduates that
covered a bit of Category Theory in the end --- see Section
\ref{teaching-adjunctions} and \cite{OchsWLD2019} --- and this forced
me to invent new techniques for working in two different styles in
parallel: a style ``for adults'', more general, abstract, and formal,
and another ``for children'', with more diagrams and examples. After
some semesters, and after writing most of the material that became
\cite{PH1}, I tried to read again some parts of Johnstone's ``Sketches
of an Elephant'', a book that always felt quite impenetrable to me,
and I found a way to present geometric morphisms in toposes to
``children''. It was based on this diagram,
\pu
$$
  \diag{gm-for-adults}
  \qquad
  \def\LG{\pshAargs{G_2}{G_3}{G_4}{G_5}}
  \def\G {\pshBargs{G_1}{G_2}{G_3}{G_4}{G_5}{G_6}}
  \def\H {\pshAargs{H_2}{H_3}{H_4}{H_5}}
  \def\RH{\pshBargs{H_2{×_{H_4}}H_3}{H_2}{H_3}{H_4}{H_5}{1}}
  \scalebox{0.6}{$
  \diag{gm-for-children}
  $}
$$
that we will discuss in detail in \ref{gms-for-children}. Its left
half is a generic geometric morphism (``for adults''), and its right
half is a very specific geometric morphism (``for children'') in which
everything is easy to understand and to visualize, and that turns out
to be ``archetypal enough''.

I showed that to the few categorists with whom I had contact and the
feedback that I got was quite positive. A few of them --- the ones who
were strictly ``adults'' --- couldn't understand why I was playing
with particular cases, and even worse, with finite categories, instead
of proving things in the most general case possible, but some others
said that these ideas were very nice, that they knew a few bits about
geometric morphisms but those bits didn't connect well, and that now
they had a family of particular cases to think about, and they had
much more intuition than before.


That was the first time that my way of using diagrams yielded
something so nice! This was the excuse that I needed to organize a
workshop on diagrammatic languages and ways to use particular
cases; here's how I advertised it (from \cite{OchsLucatelli}):
\begin{quotation}

  When we explain a theorem to children --- in the strict sense of the
  term --- we focus on concrete examples, and we avoid
  generalizations, abstract structures and infinite objects.

  When we present something to ``children'', in a wider sense of the
  term that means ``people without mathematical maturity'', or even
  ``people without expertise in a certain area'', we usually do
  something similar: we start from a few motivating examples, and then
  we generalize.

  One of the aims of this workshop is to discuss techniques for {\sl
    particularization} and {\sl generalization}. Particularization is
  easy; substituing variables in a general statement is often enough
  to do the job. Generalization is much harder, and one way to
  visualize how it works is to regard particularization as a
  projection: a coil projects a circle-like shadow on the ground, and
  we can ask for ways to ``lift'' pieces of that circle to the coil
  continously. {\sl Projections} lose dimensions and may collapse
  things that were originally different; {\sl liftings} try to
  reconstruct the missing information in a sensible way. There may be
  several different liftings for a certain part of the circle, or
  none. Finding good generalizations is somehow like finding good
  liftings.

  The second of our aims is to discuss {\sl diagrams}. For example, in
  Category Theory statements, definitions and proofs can be often
  expressed as diagrams, and if we start with a general diagram and
  particularize it we get a second diagram with the same shape as the
  first one, and that second diagram can be used as a version ``for
  children'' of the general statement and proof. Diagrams were for a
  long time considered second-class entities in CT literature
  (\cite{Kromer} discusses some of the reasons), and were omitted;
  readers who think very visually would feel that part of the work
  involved in understanding CT papers and books would be to
  reconstruct the ``missing'' diagrams from algebraic statements.
  Particular cases, even when they were the motivation for the general
  definition, are also treated as somewhat second-class --- and this
  inspires a possible meaning for what can call ``Category Theory for
  Children'': to start from the diagrams for particular cases, and
  then ``lift'' them to the general case. Note that this can be done
  outside Category Theory too; \cite{Jamnik} is a good example.

  Our third aim is to discuss {\sl models}. A standard example is that
  every topological space is a Heyting Algebra, and so a model for
  Intuitionistic Predicate Logic, and this lets us explain visually
  some features of IPL. Something similar can be done for some modal
  and paraconsistent logics; we believe that the figures for that
  should be considered more important, and be more well-known.

\end{quotation}


This is from the second announcement:

\begin{quotation}

  If we say that categorical definitions are ``for adults'' - because
  they may be very abstract - and that particular cases, diagrams, and
  analogies are ``for children'', then our intent with this workshop
  becomes easy to state. ``Children'' are willing to use ``tools for
  children'' to do mathematics, even if they will have to translate
  everything to a language ``for adults'' to make their results
  dependable and publishable, and even if the bridge between their
  tools ``for children'' and ``for adults'' is somewhat defective,
  i.e., if the translation only works on simple cases...

  We are interested in that {\sl bridge} between maths ``for adults''
  and ``for children'' in several areas. Maths ``for children'' are
  hard to publish, even informally as notes (see this thread

  \msk

  \centerline{\url{http://angg.twu.net/categories-2017may02.html}}

  \msk

  \noindent in the Categories mailing list), so often techniques are
  rediscovered over and over, but kept restricted to the ``oral
  culture'' of the area.

  Our main intents with this workshop are:

  \begin{itemize}

     \item to discuss (over coffe breaks!) the techniques of the
       ``bridge'' that we currently use in seemingly ad-hoc ways,

     \item to systematize and ``mechanize'' these techniques to make
       them quicker to apply,

     \item to find ways to publish those techniques --- in journals or
       elsewhere,

     \item to connect people in several areas working in related
       ideas, and to create repositories of online resources.

  \end{itemize}

\end{quotation}

In the UniLog 2018 I was able to chat with several categorists, and
they told me that the oral culture of CT was not as I was expected: if
there are standard ways to draw diagrams they are not widely known. I
also spent two evenings with Peter Arndt working on a certain
factorization of geometric morphisms ``for children'' --- and this
made me feel that at some point I would be able to present
applications of this diagrammatic language in ``top-tier'' conferences
that would not accept works with holes.



The following quote is from the abstract of my submission (\cite{MDE})
to the ACT2019:
\begin{quotation}

  Imagine two category theorists, Aleks and Bob, who both think very
  visually and who have exactly the same background. One day Aleks
  discovers a theorem, $T_1$, and sends an e-mail, $E_1$, to Bob,
  stating and proving $T_1$ in a purely algebraic way; then Bob is
  able to reconstruct by himself Aleks's diagrams for $T_1$ exactly as
  Aleks has thought them. We say that Bob has reconstructed the
  {\it missing diagrams} in Aleks's e-mail.

  Now suppose that Carol has published a paper, $P_2$, with a theorem
  $T_2$. Aleks and Bob both read her paper independently, and both
  pretend that she thinks diagrammatically in the same way as them.
  They both ``reconstruct the missing diagrams'' in $P_2$ in the same
  way, even though Carol has never used those diagrams herself.

\end{quotation}
%
%
and this from my submission (\cite{OchsTallinnAbs}) to Diagrams 2020:
\begin{quotation}

  Category Theory gives the impression of being an area where most
  concepts and arguments are stated and formalized via diagrams, but
  this is not exactly true... in most texts almost everything is done
  algebraically, and the reader is expected to be able to reconstruct
  the ``missing diagrams'' by himself.

  I used to believe, as an outsider, that some people who grew up
  immersed the oral culture of the area would know several techniques
  for ``drawing the missing diagrams''. My main intent when I
  organized the workshop ``Logic for Children'' at the UniLog 2018
  \cite{OchsLucatelli} was to collect some of these folklore
  techniques, compare them with the ones that I had developed myself
  to study CT, and formalize them all --- but what I found instead was
  that everybody that I could get in touch with used their own ad-hoc
  techniques, and that what I was trying to do was either totally new
  to them, or at least new in its level of detail.

\end{quotation}

The story will continue at the end of sec.\ref{the-conventions}, after
the list of conventions.

\section{The conventions \DONE}
\label{the-conventions}

The conventions that I will present now are the ones that we will need
to interpret the diagram below, that is essentially the Proposition 1
in the proof of the Yoneda Lemma in \cite[Section III.2]{CWM2}; we
will call that diagram the ``Basic Example'' and also ``diagram
$\Yzero$''. In the sections \ref{extensions}--\ref{gms-for-children}
we will see how extend DL to make it able to interpret the diagram for
geometric morphisms of the Introduction.
%
%
\pu
$$\scalebox{2.0}{$
  \slowdiag{Basic-Example}
  $}
$$

\begin{itemize}

\item[(CD)] Our diagrams are made of components that are nodes and
  arrows. The nodes can contain arbitrary expressions. The arrows work
  as connectives, and each arrow can be interpreted as the top-level
  connective in the smallest subexpression that contains it. For
  example, the curved arrow in the diagram above can be interpreted
  as:
  $$(A\ton{η}RC)↔(\catB(C,-) \ton{T} \catA(A,R-)).
  $$

\item[(C$→$)] Arrows that look like `$→$' (\qqco{\\to}) represent
  hom-sets, or, in $\Set$, spaces of functions. When a `$→$' arrow is
  named the name stands for an element of that hom-set. For example,
  in $A \ton{η} RC$ we have $η:A→RC$.

\item[(C$↦$)] Arrows that look like `$↦$' (\qqco{\\mapsto}) represent
  internal views of functions or functors. This has some subtleties;
  see Section \ref{internal-views}.

\item[(C$↔$)] Arrows that look like `$↔$' (\qqco{\\leftrightarrow})
  represent bijections or isomorphisms.

\item[(CAI)] ``Above'' usually means ``inside'', or ``internal view''.
  In the diagram above the morphism $η:A→RC$ is in $\catA$ and $C$ is
  an object of $\catB$. Also, the arrow $C \mapsto RC$ is above $\catB
  \ton{R} \catA$, and this means that it is an internal view of the
  functor $R$. Note that {\sl usually} is not {\sl always} --- and
  $\catB \ton{R} \catA$ is not an internal view of $\catB(C,-) \ton{T}
  \catA(A,R-)$.

\item[(CO)] When the definition of a component of our diagram is
  ``obvious'' in the sense of ``there is a unique natural construction
  for an object with that name'', we will usually omit its definition
  and {\sl pretend} that it is obvious; same for its uniqueness. See
  Section \ref{to-deserve-a-name}.

\item[(CC)] Everything commutes by default, and non-commutative cells
  have to be indicated explicitly. See Section \ref{freyd-notation}.

\item[(CTL)] The default ``meaning'' for a diagram {\sl without
    quantifiers} is the definition of its top-level component. There
  is a natural partial order on the components of a diagram, in which
  $α \prec β$ iff $α$ is ``more basic'' than $β$, or, in other words,
  if $α$ needs to be defined before $β$. In the diagram above the
  top-level component is the curved bijection.

\item[(CMQ)] The default ``meaning'' for a diagram with quantifiers is
  a proposition. See Sections
  \ref{freyd-notation}--\ref{freyd-with-functors} for how to obtain
  that proposition.

\item[(CAdj)] {\sl I use shapes based on my way of drawing adjunctions
  whenever possible.} I like adjunctions so much that when I want to
  explain Category Theory to someone who knows just a little bit of
  Maths I always start by the adjunction $({×}B)⊣(B{→})$ of Section
  \ref{internal-view-adjunction}; I always draw it in a canonical way,
  with the left adjoint going left, the right adjoint going right, and
  the morphisms going down. In Proposition 1 of \cite[Section
    III.2]{CWM2} the map $η$ is a universal arrow, and someone who
  learns adjunctions first sees the unit maps $η:A→(B{→}(A{×}B))$ as
  the first examples of universal arrows --- so that's why the upper
  part of the diagram above is drawn in this position.

\item[(CYo)] {\sl I use shapes based on my way of drawing the Yoneda
    Lemma whenever possible.} Look at the sections
  \ref{basic-example-as-skel}--\ref{basic-example-full} and
  \ref{yoneda-lemma}--\ref{representable-functors}: the upper parts of
  their diagrams look like parts of adjunctions, but the other parts
  do not. For example, I draw ``The functor $U:\Ring→\Set$ is
  representable'' as:
%
%
\pu
\sa{Y3+diag-Z[x]}{
                           \sa{A1}{1}
  \sa{A2}{\Z[x]}           \sa{A3}{U(\Z[x])}
  \sa{B0}{\Ring}           \sa{B1}{\Set}
  \sa{C0}{\Ring(\Z[x],-)}
                           \sa{C3}{U}
  \sa{A13}{\sm{\nameof{x}\\\univ}}
  \sa{B01}{U}
  \sa{C03}{}
  \diag{Y3+diag}
  }
$$\ga{Y3+diag-Z[x]}
$$

\item[(CDT)] A diagram acts a dictionary of default types for
  symbols. See sec.\ref{omitting-types}.

\item[(CIA)] Default types allow us to use indefinite articles in a
  precise way. An example: we have $η:\Hom_\catA(A,RC)$, so ``an $η$''
  means ``an element of $\Hom_\catA(A,RC)$''. See
  sec.\ref{indefinite-articles}.

\item[(COT)] We use a notation as close to the original text as
  possible, especially when we are trying to draw the missing diagrams
  for some existing text. If we were drawing the missing diagrams for
  the Proposition 1 of \cite[Section III.2]{CWM2} our diagram would be
  this:
  %
%
$$\pu
  \slowdiag{yoneda-cwm-0-small}
$$
but I hate Mac Lane's choice of letters, so I decided to use another
notation here.

\item[(CSk)] Suppose that we have a piece of text --- say, a paragraph
  $P$ --- and we want to reconstruct the ``missing diagram'' $D$ for
  $P$. Ideally this $D$ should be a ``skeleton'' for $P$, in the sense
  that it should be possible to reconstruct the ideas in $P$ from the
  diagram $D$ using very few extra hints; see \cite[sec.12]{IDARCT}.

\item[(CTT)] Our diagrams should be close to Type Theory: it should be
  possible to use them as a scaffolding for formalizing our text in
  (pseudocode for) a proof assistant.

\item[(CFSh)] The image by a functor of a diagram $D$ is drawn with
  the same shape as $D$.

\item[(CISh)] The internal view of a diagram $D$ is drawn with the
  same shape as $D$, modulo duplications --- see section
  \ref{internal-views}.

\item[(CPSh)] A particular case of a diagram $D$ is drawn with the
  same shape as $D$.

\item[(CNSh)] A translation of a diagram $D$ to another notation is
  drawn with the same shape as $D$.

\end{itemize}

The conventions (CD)--(CMQ) and (CFSh)--(CNSh) all appear in diagrams
in \cite{MacLaneNotes}, \cite{Freyd76}, \cite{FreydScedrov},
\cite{Taylor}, \cite{Riehl}, \cite{Leinster}, but very few of them are
spelled out explicitly, and the idea of ``same shape'' is never
stressed. See \cite[p.179]{NederpeltGeuvers} for a neat example of
``substitution produces something with the same shape'' and
\cite{PenroseSIGGRAPH2020} for a language for drawing diagrams from
high-level specifications in which it may be possible to implement the
rules about ``same shape''.

The other conventions {\sl may} be new, but remember from the
introduction that most of the work on diagrammatic languages lies
below the threshold of publishability... so conventions corresponding
to those may be folklore knowledge in groups that I don't have contact
with {\sl yet}.

The convention (COT) is obvious in retrospect, but giving a name to it
saved me from my tendency to invent new notations. The conventions
(CDT) and (CIA) replace the idea of downcasings from
\cite[sec.3]{IDARCT}, that didn't work well. Sections
\ref{extensions}--\ref{gms-for-children} show how to add new
conventions to DL, and sec.\ref{opposite-categories} shows that we can
add a bad convention and mark it as a temporary hack.

There are many notations for Type Theory. To make this chapter more
readable in the convention (CTT) I will use a pseudocode that is
halfway between standard mathematical notation and Agda; the companion
paper \cite{MissingAgda} will show how to translate it to real Agda
(with the library \cite{HuCarette}).

Most texts on CT use diagrams to {\sl prove} theorems. Here will use
them to {\sl understand} theorems, and to translate between languages.
Our approach can be seen as an extension of \cite{Ganesalingam} to
Category Theory; see also \cite{DSLsofMath}, that is a recent book
that follows many of the ideas in \cite{Ganesalingam}.


%
%

%
%
\section{Finding ``the'' object with a given name \DONE}
\label{to-deserve-a-name}

One of the books that I tried to read when I was starting to learn
Category Theory was Mac Lane's \cite{CWM2}. It is written for readers
who know a lot of mathematics and who can follow some steps that it
treats as obvious. I was not (yet) a reader like that, but I wanted to
become one.

There is one specific thing that \cite{CWM2} does pretending that it
is obvious that I found especially fascinating. It ``defines''
functors by describing their actions on objects, and it leaves to the
reader the task of discovering their actions on morphisms. Let's see
how to find these actions on morphisms.

A functor $F:\catA→\catB$ has four components:
$$F=(F_0, F_1, \respids_F, \respcomp_F).$$
They are its action on objects, its action on morphisms, the assurance
that it takes identity maps to identity maps, and the assurance that
it respects compositions. When Mac Lane says this,
\begin{quote}
Fix a set $B$. Let $(×B)$ denote {\sl the} functor that takes each set
$A$ to $A×B$.
\end{quote}
he is saying that $(×B)_0 A = A×B$, or, more precisely, this:
$$(×B)_0 := λA.\,A×B$$

The ``{\sl the}'' in the expression ``Let $(×B)$ denote {\sl the}
functor...'' implies that the precise meaning of $(×B)_1$ is easy to
find, and that it is easy to prove $\respids_{(×B)}$ and
$\respcomp_{(×B)}$.

If $f:A'→A$ then $(×B)_1 f : (×B)_0 A' → (×B)_0 A$. We know the {\sl
  name} of the image morphism, $(×B)_1 f$, and its {\sl type},
$$(×B)_1 f : A'×B → A×B,$$
and it is implicit that there is an ``obvious'' natural construction
for this $(×B)_1 f$ from $f$. A natural construction is ---
TA-DAAAA!!! --- a $λ$-term, so we are looking for a term of type $A'×B
→ A×B$ that can be constructed from $f:A'→A$.

In a big diagram:
\pu
$$\ded{foo1} \quad ⇒ \ded{foo2}$$

A double bar in a derivation means ``there are several omitted steps
here'', and sometimes a double bar suggests that these omitted steps
are obvious. The derivation on the left says that there is an
``obvious'' way to build a $(×B)_1f:A'{×}B→A{×}B$ from a
``hypothesis'' $f:A'→A$. If we expand its double bar we get the tree
at the right, that shows that the ``precise meaning'' for $(×B)_1f$ is
$(λp⠆A'{×}B.(f(πp),π'p)$. More formally (and erasing a typing),
$$(×B)_1 := λf.(λp.(f(πp),π'p)).$$

The expansion of the double bar above becomes something more familiar
if we translate the trees to Logic using Curry-Howard:
\pu
$$\ded{foo-logic1} \quad ⇒ \ded{foo-logic2}$$

We obtain the tree at the right by {\sl proof search}.

Let's give a name for the operation above that obtained a term of type
$A'×B→A×B$: we will call that operation {\sl term search}, or, as it
is somewhat related to type inference, {\sl term inference}.

Term search may yield several different construction and trees, and so
several non-equivalent terms of the desired type. When Mac Lane says
``{\sl the} functor $(×B)$'' he is indicating that:

\begin{itemize}

\item a term for $(×B)_1$ is easy to find (note that we use the
  expression ``a {\sl precise meaning} for $(×B)_1$''),

\item all other natural constructions for something that ``deserves
  the name'' $(×B)_1$ yield terms equivalent to that first, most
  obvious one,

\item proving $\respids_{(×B)}$ and $\respcomp_{(×B)}$ is trivial.

\end{itemize}

In many situations we will start by just the name of a functor, as the
``$(×B)$'' in the example above, and from that name it will be easy to
find {\sl the} ``precise meaning'' for $(×B)_0$, and from that the
``precise meaning'' for $(×B)_1$, and after that proofs that
$\respids_{(×B)}$ and $\respcomp_{(×B)}$. We will use the expression
``...deserving the name...'' in this process --- terms for $(×B)_0$,
$(×B)_1$, $\respids_{(×B)}$, and $\respcomp_{(×B)}$ ``deserve their
names'' if they obey the expected constraints.

For a more thorough discussion see \cite{IDARCT}.

\msk

These ideas of ``finding a precise meaning'' and ``finding (something)
deserving that name'' can also be applied to morphisms, natural
transformations, isomorphisms, and so on.

In Section \ref{basic-example-NTs} we will see how to find natural
constructions for the two directions of the bijection in the Basic
Example --- or how the expand the double bars in the two derivations
here:
$$\pu
  \slowdiag{Basic-Example}
  \qquad
  \begin{array}{c}
  \ded{bij-down}
  \\ \\
  \ded{bij-up}
  \end{array}
$$

\msk

I am not aware of any papers or books on CT that discuss how to
(re)construct a functor from its action on objects or from its name,
but Agda has a tool that can be used for that: look for the section on
``Automatic Proof Search'' in \cite{AgdaUserManual}.



%

\section{Freyd's diagrammatic language \DONE}
\label{freyd-notation}

In \cite{Freyd76} Peter Freyd presents a very nice diagrammatic
language that can be used to express {\sl some} definitions from
Category Theory. For example, this is the statement that a category
has all equalizers:
$$\pu
  \scalebox{0.8}{$
  \diag{cat-has-equalizers}
  $}
$$

All cells in these diagrams commute by default, and non-commuting
cells have to be indicated with a `?'. Each vertical bar with a `$∀$'
above it means ``for all extensions of the previous diagram to this
one such that everything commutes''; a vertical bar with a `$∃!$'
above it means ``there exists a unique extension of the previous
diagram to this one such that everything commutes'', and so on. See
the scan in \cite{Freyd76} for the basic details of how to formalize
these diagrams, and the book \cite[p.28 onwards]{FreydScedrov}, for
tons of extra details, examples, and applications.

Let's call the subdiagrams of a diagram like the one above its
``stages''. Its stage 0 is empty, its stage 1 has two objects and two
arrows, its last stage has four objects and five arrows, and the
quantifiers separating the stages are $Q_1=∀$, $Q_2=∃$, $Q_3=∀$,
$Q_4=∃!$. They are structured like this:
%
\pu
$$
  \diag{freyd-stages-1}
$$

I was not very good at drawing all stages separately --- it was
boring, it took me too long, and I often got distracted and committed
errors --- so I started to play with extensions of that diagrammatic
language.

%
\subsection{Adding quantifiers \DONE}
\label{freyd-with-quantifiers}

Here is a simple way to draw all stages at once. We start from a
diagram for the ``last stage with quantifiers'', that we will call
$LSQ$:
%
$$\pu
  \scalebox{1.75}{$
  \diag{cat-has-equalizers-with-quants}
  $}
$$

We can recover all the stages and quantifiers from it. The numbered
quantifiers in it are $∀_1$, $∃_2$, $∀_3$, and $∃!_4$. The highest
number in them is 4, so we set $n=4$ ($n$ is the index of the last
stage), and we set ``stage 4 with quantifiers'', $SQ_4$, to $LSQ$. To
obtain the $SQ_3$ from $SQ_4$ we delete all nodes an arrows in $SQ_4$
that are annotated with a `$∃!_4$'; to obtain $SQ_2$ from $SQ_3$ we
delete all nodes an arrows in $SQ_3$ that are annotated with a
`$∀_3$', and so on until we get a diagram $SQ_0$, that in this example
is empty. To obtain each $S_k$ --- a stage in the original
diagrammatic language from Freyd, that doesn't have quantifiers ---
from the corresponding $SQ_k$ we treat all the quantifiers in $SQ_k$
as mere annotations, and we erase them; for example, `$∃_2e$' becomes
`$e$', and $∀_1A$ becomes $A$. To obtain the quantifiers $Q_1$, $Q_2$,
$Q_3$, $Q_4$ that are put in the vartical bars that separate the
stages, we just assign $∀_1$, $∃_2$, $∀_3$, and $∃!_4$ to them,
without the numbers in the subscripts.

Bonus convention: when the quantifiers in a diagram are just `$∀$'s
and `$∃!$'s without subscripts the `$∀$'s are to be interpreted as
`$∀_1$' and the `$∃!$'s as `$∃!_2$'s.

%
\subsection{Adding functors \DONE}
\label{freyd-with-functors}

Freyd's language can't represent functors --- in the sense of diagrams
like the ones in sec.\ref{internal-view-functor} --- and I wanted to
use it to draw the missing diagrams for definitions involving
functors, so I had to extend it again.


Let me use an example to discuss this. This is the definition of
universal arrow in \cite[p.55]{CWM2}, including the original diagram,
modulo change of letters:

%

\begin{quotation}

  {\bf Definition.} If $R: \catB→\catA$ is a functor and $A$ an object
  of $\catA$, a universal arrow from $A$ to $R$ is a pair $(B,η)$
  consisting of an object $B$ of $\catB$ and and arrow $η:A→RB$ of
  $\catA$ such that to every pair $(B',g)$ with $B'$ an object of
  $\catB$ and $g:A→RB'$ an arrow of $\catA$, there is a unique arrow
  $f:B→B'$ of $\catB$ with $Rf∘η=g$. In other words, every arrow $h$
  to $R$ factors uniquely through the universal arrow $η$, as in the
  commutative diagram:
  $$\pu
    \diag{univ-arrow-cwm-my-letters}
  $$

\end{quotation}

The definition itself goes only up to the ``with $Rf∘η=g$.'', so let
me ignore the part starting from ``In other words'', and draw a better
``missing diagram'' for the definition:
%
$$\pu
  \diag{universal-arrow-stages}
$$

This diagram is quite close to being a skeleton for the definition of
universal arrow. It can be interpreted as a proposition, and the only
extra hint that we need is that ``universalness'' for the arrow $η$
corresponds to the truth of that proposition. Here's how to extract
the proposition from it:
$$\begin{array}{rl}
  \text{In a context where:}
    & \catA \text{ is a category}, \\
    & \catB \text{ is a category}, \\
    & R:\catB \to \catA, \\
    & A ∈ \catA, \\
    & B ∈ \catB, \\
    & η:A→RB, \\
  \text{for all}
    & B'∈\catB \text{ and} \\
    & g:A→RB', \\
  \text{there exists a unique}
    & f:B→B' \text { such that} \\
    & Rf∘η=g. \\
  \end{array}
$$

To convert that to a definition of universalness we just have to
replace the ``for all'' by ``$(B,η)$ is a universal arrow for $A$ to
$R$ iff for all''.

The convention for quantifiers from sec.\ref{freyd-with-quantifiers}
lets us rewrite the diagram in three stages above as:
%
$$\pu
  \scalebox{1.5}{$
  \diag{universal-arrow-quants}
  $}
$$

Also, I noticed that I could omit most typings when they could be
inferred from the diagram. I could ``formalize'' the diagram above as:
``in a context where $(\catA, \catB, R, A, B, η)$ are as in the
diagram above, we say that $(B,η)$ is a universal arrow from $A$ to
$R$ when $∀(B',g).∃!f.(Rf∘η=g)$''. This looked too loaded to be used
in public, but it was practical for private notes --- and I could even
omit the ``$Rf∘η=g$'', as everything commutes by default. In
sec.\ref{omitting-types} we will see a way to formalize this method
for omitting and reconstructing types, and in
sec.\ref{indefinite-articles} we will see a second way to define
universalness.

\bsk

Note that when we erase a node or arrow we also erase everything that
depends on it. In the example above $SQ_2$ has an arrow labeled
$∃!_2f$; to obtain $SQ_1$ from $SQ_2$ we have to erase that arrow, the
arrow $Rf$, and the arrow $f \mapsto Rf$ --- and to obtain $SQ_0$ from
$SQ_1$ we have to erase the arrow $g$, the node $B'$, the node $RB'$,
and the arrow $B' \mapsto RB'$.

%
\section{Internal views \DONE}
\label{internal-views}

My favorite way of introducing internal views is with the diagram
below:
%
%
\def\ooo(#1,#2){\begin{picture}(0,0)\put(0,0){\oval(#1,#2)}\end{picture}}
\def\oooo(#1,#2){{\setlength{\unitlength}{1ex}\ooo(#1,#2)}}
%
\pu
$$\begin{array}{rrcl}
   \sqrt{\;\;}: & \N &→& \R \\
                &  n &↦& \sqrt{n} \\
   \end{array}
   \qquad
   \diag{second-blob-function}
$$

\def\longmapsto{\diagxyto/|->/}

The parts with the two blobs and `$\longmapsto$'s between them is
based on how I was taught sets and functions when I was a kid; it is
an internal view of the $\N \ton{\sqrt{\phantom{a}}} \R$ below it. Not
all elements of $\N$ are shown in the blob-view of $\N$, but the ones
that are shown are named; compare this with \cite[p.2
  onwards]{LawvereRosebrugh}, in which the elements are usually dots.

The arrow $n \longmapsto \sqrt{n}$ between the blobs shows a {\sl
  generic element} of $\N$ and its image, and the other
`$\longmapsto$'s are {\sl substitution instances of it}, like this:
$$(n \longmapsto \sqrt{n}) [n:=2] = (2 \longmapsto \sqrt{2})$$

In some cases, like $4 \longmapsto 2$, we write 2 instead of
$\sqrt{4}$ because $\sqrt{4}$ ``reduces to'' 2, as explained in the
next section.

%
%
\subsection{Reductions \DONE}
\label{reductions}

\def\squigton#1{\overset{#1}{\squigto}}

The convention (C$\mapsto$) says that an arrow $α \mapsto β$ above an
arrow $A \ton{f} B$ should be interpreted as meaning $f(α) \squigto
β$, where `$\squigto$' means ``reduces to''; the standard example is
$\sqrt{4} \squigto 2$. In a diagram:
%
$$\pu
  \diag{reductions-0}
$$

The idea of reduction comes from $λ$-calculus. We write $α
\squigton{1} β$ to say that the term $α$ reduces to $β$ in one step,
and $α \squigton{*} γ$ to say that there is a finite sequence of
one-step reductions that reduce $α$ to $γ$. Here we are interested in
reduction in a system with constants, in which for example
$(\sqrt{\phantom{a}})(4) \squigton{1} 2$.

Here is a directed graph that shows all the one-step reductions
starting from $g(2+3)$, considering $g(a) = a·a+4$:
%
%
$$\pu
  \diag{reduce-g}
$$

Note that all reductions sequences starting from $g(2+3)$ terminate at
the same term, 29 --- ``the term $g(2+3)$ is strongly normalizing''
---, and reduction sequences from $g(2+3)$ may ``diverge'' but they
``converge'' later --- this is the ``Church-Rosser Property'', a.k.a.
``confluence''.

A good place to learn about reduction in systems with constants is
\cite{SICP}.



%
%
\subsection{Functors \DONE}
\label{internal-view-functor}


By the convention (CFSh) the image of the diagram above $\catA$ in the
diagram below --- remember that {\sl above} usually means {\sl inside}
---
%
$$\pu
  \diag{internal-view-functor-0}
$$
is a diagram with the same shape over $\catB$. We draw it like this:
%
$$\pu
  \diag{internal-view-functor-1}
$$

In this case we don't draw the arrows like $A_1 \mapsto FA_1$ because
there would be too many of them --- we leave them implicit.

We say that the diagram above is {\sl an} internal view of the functor
$F$. To draw {\sl the} internal view of the functor $F: \catA → \catB$
we start with a diagram in $\catA$ that is made of two generic objects
and a generic morphism between them. We get this:

$$\pu
  \diag{internal-view-functor-2}
$$

Compare this with the diagram with blob-sets in Section
\ref{internal-views}, in which the `$n \mapsto \sqrt{n}$' says where a
generic element is taken.

Any arrow of the form $α \mapsto β$ above a functor arrow $\catA
\ton{F} \catB$ is interpreted as saying that $F$ takes $α$ to $β$, or,
in the terminology of the section \ref{reductions}, that $Fα$ reduces
to $β$. So this diagram 
%
$$\pu
  \diag{internal-view-functor-3}
$$
defines $(A×)$ as:
$$\begin{array}{rcl}
  (A×)_0 &:=& λB.\,A×B,\\
  (A×)_1 &:=& λf.λp.(πp,f(π'p)).\\
  \end{array}
$$

In this case we can also use internal views of $(A×)$ to define
$(A×)_1$:
%
$$\pu
  \diag{internal-view-functor-4}
$$

%
%
\subsection{Natural transformations \DONE}
\label{internal-view-NT}

Suppose that we have two functors $F,G:\catA → \catB$ and a natural
transformation $T:F→G$. A first way to draw an internal view of $T$ is
this:
%
$$\pu
  \diag{internal-view-NT-0}
$$

If we start with a morphism $h:C→D$ in $\catA$, like this,
%
%
$$\pu
  \diag{NT-00}
$$
the convention (CFSh) would yield an image of $h$ by $F$ and another
by $G$, and we can draw the arrows $TC$ and $TD$ to obtain a commuting
square in $\catB$:
%
%
$$\pu
  \diag{NT-0}
$$

This way of drawing internal views of natural transformations yields
diagrams that are too heavy, so we will usually draw them as just
this:
%
%
$$\pu
  \diag{NT-1}
$$
Note that the input morphism is at the left, and above $F \ton{T} G$
we draw its images by $F$, $G$, and $T$.

When the codomain of $F$ and $G$ is $\Set$ we will sometimes also draw
at the right an internal view of the commuting square, like this:
%
$$\pu
  \diag{NT-2}
$$
Then the commutativity of the middle square is equivalent to
$∀x∈FC.(Gh∘TC)(x)=(TD∘Ff)(x)$. Note that in this case the square at
the right is an internal view of an internal view.

In Section \ref{to-deserve-a-name} we saw that a functor has four
components. A natural transformation has two: $T=(T_0, \sqcond_T)$,
where $T_0$ is the operation $C \mapsto TC$ and $\sqcond_T$ is the
guarantee that all the induced squares commute.

%


%
%
\subsection{Adjunctions \DONE}
\label{internal-view-adjunction}


We will draw adjunctions like this,
%
$$\pu
  \diag{generic-adjunction}
$$
with the left adjoint going left and the right adjoint going right. My
favorite names for the left and right adjoints are $L$ and $R$. The
standard notation for that adjunction is $L⊣R$.

The top-level component of the diagram above is the bijection arrow in
the middle of the square --- it says that $\Hom(LA,B) ↔ \Hom(A,RB)$.
It is implicit that we have bijections like that for all $A$ and $B$;
it is also implicit that that bijection is natural in some sense.

We will sometimes expand adjunction diagrams by adding unit and counit
maps, the unit and the unit as natural transformations, the actions of
$L$ and $R$ on morphisms, and other things. For example:
%
$$\pu
  \diag{generic-adjunction-with-with}
$$

\def\HomA#1{HomA(#1)}
\def\HomB#1{HomB(#1)}
\def\HomA#1{\catA(#1)}
\def\HomB#1{\catB(#1)}

We can obtain the naturality conditions by regarding $♭$ and $♯$ as
natural transformations and drawing the internal views of their
internal views:
%
$$\pu
  \diag{adj-nat-conditions}
$$

%

\subsection{A way to teach adjunctions \DONE}
\label{teaching-adjunctions}

I mentioned in the first sections that I have tested some parts of
this language in a seminar course --- described here:
\cite{OchsWLD2019} --- and that in it I taught Categories starting by
adjunctions. Here's how: we started by the basics of $λ$-calculus and
some sections of \cite{PH1}, and then I asked the students to define
each one of the operations in the right half of the diagram below as
$λ$-terms:
$$
  \pu
  \diag{generic-adjunction-big}
  \qquad
  \diag{(xB)-|(B->)}
$$

Then we saw the definition of functors, natural transformations and
adjunctions, and we checked that the right half is a particular case
(``for children''!) of the diagram for a generic adjunction in the
left half. After that, and after also checking that in the Planar
Heyting Algebras of \cite{PH1} we have an adjunction $(∧Q)⊣(Q→)$, I
helped the students to decypher some excerpts of \cite{Awodey}.

\msk

From the components of the generic adjunction in the diagram above it
is possible to build this big diagram:
$$\pu
  \diag{teaching-adjunctions-1}
$$

Let's use these names for its subdiagrams: $\sm{ A \\ BCDEF \\ G \\ I}$.


A {\sl fully-specified adjunction} between categories $\catB$ and
$\catA$ has lots of components: $(L, R, ε, η, ♭, ♯, \univ(ε),
\univ(η))$, and maybe even others, but usually we define only some of
these components; there is a Big Theorem About Adjunctions (below!)
that says how to reconstruct the fully-specified adjunction from some
of its components.

Some parts of the diagram above can be interpreted as definitions,
like these:
$$\begin{array}{c}
  Lf := (η_A∘f)^♭ \\
  [5pt]
  g := ε_B∘Lh 
    \qquad ε_B := (\id_{RB})^♭
    \qquad η_A := (\id_{LA})^♯
    \qquad h := Rg∘η_A \\
  [5pt]
  Rk := (k∘η_B)^♯ \\
  \end{array}
$$

The subdiagrams $B$ and $F$ can also be interpreted in the opposite
direction, as:
$$\begin{array}{rclcrcl}
  g^♯ &:=& (∀A.∀g.∃!h)Ag    &\phantom{mmmmmm}&  h^♭ &:=& (∀B.∀h.∃!g)Bh \\
       &=& (\univ_{ε_B})Ag  &&                       &=& (\univ_{η_A}) Bh \\
  \end{array}
$$

The notations $(∀A.∀g.∃!h)Ag$ and $(\univ_{ε_B})Ag$ are clearly abuses
of language --- but it's not hard to translate them to something
formal, and these notations inspired great discussions in the
classroom... also, they can help us to understand and formalize
constructions like this one,
%
$$\pu
  Lf := (\univ_{η_A})(LA)(η_A∘f)
  \qquad
  \diag{building-L_1f}
$$
that are needed in cases like the part (ii) of the Big Theorem.

The Big Theorem About Adjunctions is this --- it's the Theorem 2 in
\cite[page 83]{CWM2}, but with letters changed to match the ones we
are using in our diagrams:

\def\ORIG#1{\msk\ColorBrown{#1}}
\def\ORIG#1{}

\newpage

\begin{quotation}

  \ORIG{{\bf Theorem 2.} Each adjunction $〈F,G,φ〉: X \rightharpoonup
    A$ is completely determined by the items in any one of the
    following lists:}

  {\bf Big Theorem About Adjunctions.} Each adjunction $〈L,R,♯〉: \catA
  \rightharpoonup \catB$ is completely determined by the items in any
  one of the following lists:

  \ORIG{(i) Functors $F$, $G$, and a natural transformation $η: 1_X
    \tnto GF$ such that each $η_x: x→GFx$ is universal to $G$ from
    $x$. Then $φ$ is defined by (6).}

  (i) Functors $L$, $R$, and a natural transformation $η:
  \id_\catA→RL$ such that each $η_A: A→RLA$ is universal to $R$ from
  $A$. Then $♯$ is defined by (6).

  \ORIG{(ii) The functor $G: A → X$ and for each $x∈X$ an
    object $F_0x∈A$ and a universal arrow $η_x:x→GF_0x$ from $x$
    to $G$. Then the functor $F$ has object function $F_0$ and is
    defined on arrows $h:x→x'$ by $GFh∘η_x = η_{x'}∘h$.}

  (ii) The functor $R: \catB → \catA$ and for each $A∈\catA$ an object
  $L_0A∈\catB$ and a universal arrow $η_A:A→RL_0A$ from $A$ to $R$.
  Then the functor $L$ has object function $L_0$ and is defined on
  arrows $f:A'→A$ by $RLf∘η_{A'} = η_A∘f$.

  \ORIG{(iii) Functors $F$, $G$, and a natural transformation $ε: FG
    \tnto I_A$ such that each $ε_a:FGa→a$ is universal from $F$ to
    $a$. Here $φ^{-1}$ is defined by (7).}

  (iii) Functors $L$, $R$, and a natural transformation $ε:
  LR→\id_\catB$ such that each $ε_B:LRB→B$ is universal from $L$ to
  $B$. Here $♭$ is defined by (7).

  \ORIG{(iv) The functor $F:X→A$ and for each $a∈A$ an object $G_0a∈X$
    and an arrow $ε_a:FG_0a→a$ universal from $F$ to $a$.}

  (iv) The functor $L:\catA→\catB$ and for each $B∈\catB$ an object
  $R_0B∈\catA$ and an arrow $ε_B:LR_0B→B$ universal from $L$ to $B$.

  \ORIG{(v) Functors $F$, $G$ and natural transformations $η:I_x \tnto
    GF$ and $ε: FG \tnto I_A$ such that both composites (8) are the
    identity transformations. Here $φ$ is defined by (6) and $φ^{-1}$
    by (7).}

  (v) Functors $L$, $R$ and natural transformations $η:\id_\catA→RL$
  and $ε:LR→\id_\catB$ such that both composites (8) are the identity
  transformations. Here $♯$ is defined by (6) and $♭$ by (7).

\end{quotation}



\newpage

%
\section{Types for Children}
\label{types-for-children}


We will need a bit of Type Theory in the sections \ref{ness} and
\ref{comma-categories}. We will need some non-standard notational
conventions that appear more or less naturally when we present Theory
Theory ``for children'' in the right way --- let's see the details.

\msk


Section 6 of \cite{SelingerLN} has a very good presentation of types
``for adults'': it uses expressions like $A×B$ and $A→B$ as and treats
them as purely syntactical objects, but each one comes with an
``intended meaning''. Let's start by defining a version ``for
children'' of that in which these intended meanings become more
concrete, and then we will work in the version ``for children'' and in
the version ``for adults'' in parallel.

\subsection{Dependent types}
\label{dependent-types}

In our version ``for children'':

\begin{itemize}

\item all types are sets,

\item some sets are types,

\item every finite subset of $\N$ is a type,

\item if $A$ and $B$ are types then $A×B$ and $A→B$ are types. $A×B$
  is the space of pairs of the form $(a,b)$ in which $a∈A$ and $b∈B$,
  and $A→B$ is the space of functions from $A$ to $B$,

\item $a:A$ means $a∈A$ --- the distinction between `$:$' and `$∈$'
  will only appear in other settings,

\item ``space of'' means ``set of''. The space of functions from $A$
  to $B$ is the set of all functions from $A$ to $B$, and each
  function is considered as a set of input-output pairs. So, for
  example, if $A=\{2,3\}$ and $B=\{4,5\}$ then:
  $$\def\fa#1#2{\csm{(2,#1),\\(3,#2)\,}}
    \begin{array}{rcl}
    A×B &=& \{(2,4),
              (2,5),
              (3,4),
              (3,5),
            \}, \\
    A→B &=& \left\{ \fa44, \fa45, \fa54, \fa55
            \right\} \\
    \end{array}
  $$

\item if $A$ is a type and $(C_a)_{a∈A}$ is a family of types indexed
  by $A$ then $Πa{:}A.C_a$ and $Σa{:}A.C_a$ are dependent types
  defined in the usual way, and $(a{:}A) → C_a$ and $(a{:}A)×C_a$ are
  alternate notations for $Πa{:}A.C_a$ and $Σa{:}A.C_a$ (see
  \cite[section 2]{Norell08}). Formally,
  $$\def\aCa{\bigcup_{a∈A}C_a}
    \begin{array}{rcl}
     Σa{:}A.C_a &=& \setofst{(a,c)∈A × (\aCa) }{a∈A, \; c∈C_a} \\
    (a{:}A)×C_a &=& \setofst{(a,c)∈A × (\aCa) }{a∈A, \; c∈C_a} \\
     Πa{:}A.C_a &=& \setofst{f:A→(\aCa) }{∀a∈A.\;f(a)∈C_a} \\
    (a{:}A)→C_a &=& \setofst{f:A→(\aCa) }{∀a∈A.\;f(a)∈C_a} \\
    \end{array}
  $$
  If $A=\{2,3\}$, $C_2=\{6,7\}$, and $C_3=\{7,8\}$ then:
  $$\def\fa#1#2{\csm{(2,#1),\\(3,#2)\,}}
    \begin{array}{rcl}
    (a{:}A)×C_a &=& \{(2,6),
                      (2,7),
                      (3,7),
                      (3,8),
                    \}, \\
    (a{:}A)→C &=& \left\{ \fa67, \fa68, \fa77, \fa78
                  \right\}. \\
    \end{array}
  $$ 

\end{itemize}

\subsection{Witnesses}
\label{witnesses}

If $P$ is a proposition we will write $⟦P⟧$ for its {\sl space of
  witnesses}, or its {\sl space of proofs}. The exact definition of
$⟦P⟧$ will usually depend on the context, but we always have $⟦P⟧=∅$
when $P$ is false and $⟦P⟧\neq∅$ when $P$ is true. In some situations
all the witnesses of a proposition $P$ will be identified --- this is
called {\sl proof irrelevance}; see \cite[p.340]{NederpeltGeuvers} ---
and all the spaces of witnesses will be either singletons or empty
sets; in other situations some `$⟦P⟧$'s will have more than one
element.


The notation $\AProofOf{P}$ will denote a witness that $P$ is true.
Formally, $\AProofOf{P}$ is a variable (with a long name!) whose type
is $\AllProofsOf{P}$. A good way to remember this notation is that
$\AllProofsOf{P}$ looks like a box and $\AProofOf{P}$ looks like
something that comes in that box.

\msk

{

\def\a{\mathsf{a}}
\def\b{\mathsf{b}}

In Agda the operation `$≡$' returns a space of proofs of equality. If
$\a$ and $\b$ are expressions with the same type then Agda's `$\a≡\b$'
corresponds to our $\AllProofsOf{\a=\b}$, and people sometimes use the
name `$\a{≡}\b$' to denote an element of $\a≡\b$ --- we use
$\AProofOf{\a=\b}$ for that. See the section ``Equality'' in
\cite{WadlerPLFA} for simple examples, and Agda's standard library for
more examples.

}


\subsection{Judgments}


The main objects of Type Theory are {\sl derivable judgements}. A
derivable judgment is one that can appear as the root node of a
derivation in which each bar is an application of one the rules in
\cite[p.127]{NederpeltGeuvers}. These derivations are usually huge ---
for example, here is a derivation for $A⠆Θ,B⠆Θ⊢(Πa⠆A.B)⠆Θ$:
\pu
$$\resizebox{0.7\textwidth}{!}{$
    \ded{depprod1}
  $}
$$
so rarely draw them explicitly, and we use other tools to show that
certain judgments are derivable.

Every derivable judgment obeys this (taken verbatim from
\cite[p.52]{SelingerLN}):
%
%
\begin{quotation}

  A typing judgment is an expression of the form
  $$x_1:A_1, x_2:A_2, \ldots, x_n:A_n ⊢ M:A.$$
  
  Its meaning is: "under the assumption that $x_i$ is of type $A_i$ ,
  for $i=1\ldots n$, the term $M$ is a well-typed term of type $A$.'' The
  free variables of $M$ must be contained in $x_1, \ldots, x_n$.

\end{quotation}

Understanding what this means in the version ``for children'' will
take us quite close to understanding that in Type Theory ``for
adults''. We will do that in the next section.

Let me just correct a simplification. The main objects of the Type
Theory used in Agda and in most other proof assistants are derivable
judgments {\sl with definitions}, as explained in the chapters 8--10
of \cite{NederpeltGeuvers}. A judgment with definitions is written as
$Δ;Γ⊢M:N$, where $Δ$ is a list of definitions
(\cite[def.9.2.1]{NederpeltGeuvers}); we will mostly ignore the `$Δ$'
here.


\subsection{Set comprehensions}

\def\und#1#2{\underbrace{#1}_{#2}}
\def\und#1#2{\underbrace{#1}_{\text{#2}}}
\def\ug#1{#1}
\def\uf#1{#1}
\def\ue#1{#1}
\def\uc#1{#1}
\def\UG#1{\und{#1}{generator}}
\def\UF#1{\und{#1}{filter}}
\def\UE#1{\und{#1}{expression}}
\def\UC#1{\und{#1}{context}}

The part at the left of the `$⊢$' in a typing judgment is called a
{\sl typing context}. Typing contexts also appear in set
comprehensions. Let's see an example:
$$\begin{array}{rl}
    & \setofst {\ue{10a+b}}
               {\uc{\ug{a∈\{1,2\}}, \ug{b∈\{2,3\}}, \uf{a<b}}} \\[7.5pt]
  \squigto \phantom{mm}
    & \setofsc {\UC{\UG{a∈\{1,2\}}, \UG{b∈\{2,3\}}, \UF{a<b}}}
               {\UE{10a+b}} \\
  \end{array}
$$

\def\acontext{\ang{\textsf{context}}}
\def\aexpr   {\ang{\textsf{expr}}}

We rewrote the comprehension $\setofst{ \aexpr }{ \acontext }$ to
$\setofsc{ \acontext }{ \aexpr }$ for clarity, and we marked which
parts of the context act as ``generators'' and which ones act as
``filters''. The context above can be rewritten in type-theoretical
notation as:
$$   a : \{1,2\},
  \; b : \{2,3\},
  \; \AProofOf{a{<}b} : \AllProofsOf{a{<}b}
$$

A {\sl value} for that context is a triple $(a,b,\AProofOf{a{<}b})$,
where $a∈\{1,2\}$, $b∈\{2,3\}$, and $\AProofOf{a{<}b}$ is a guarantee
that $a<b$ is true.


\subsection{Omitting types}
\label{omitting-types}


The diagram at the left below is a copy of the one from
sec.\ref{freyd-with-functors}, but now we will interpret it in a
different way, as a ``dictionary of (default) types''. For example, it
says that when the symbol $η$ appears without a type its type will be
the default one given by the diagram: $η:A→RC$, or $η:\catA(A,RC)$.
The default types are listed at the right.

$$\pu
  \scalebox{1.0}{$
  \diag{universal-arrow-omitting-1}
  $}
  \qquad
  \begin{array}{rcl}
  \catA &\text{is}& \text{a category} \\
  \catB &\text{is}& \text{a category} \\
      R &:& \catA → \catB \\
      A &∈& \catA \\
      C &∈& \catB \\
      D &∈& \catB \\
      η &:& A→RC \\
      f &:& C→D \\
      g &:& A→RD \\
  \end{array}
$$

\subsection{Indefinite articles}
\label{indefinite-articles}

We will use the diagram above to redefine universalness. In our old
definition, from sec.\ref{freyd-with-functors}, universalness is just a
``property''; in our new definition it will be a pair made of a
``structure'' and a ``property'' (see sec.\ref{ness}).

Suppose that 
$$\begin{array}{rcl}
  \catA &\text{is}& \text{a category}, \\
  \catB &\text{is}& \text{a category}, \\
      R &:& \catA → \catB, \\
      A &∈& \catA, \\
      C &∈& \catB, \\
      η &:& A→RC, \\
  \end{array}
$$
and let $♯$ be this operation:
$$♯ \;\;=\;\; λD.λf.\;Rf∘η.
$$

\def\asharpis{\AProofOf   {♭ = ♯^{-1}}}
\def\psharpis{            (♭ = ♯^{-1})}
\def\Ssharpis{\AllProofsOf{♭ = ♯^{-1}}}

Note that the types of $D$ and $f$ are given by the diagram: $D$ is an
object of $\catB$, and $f:C→D$. Then ``universalness'' for the tuple
$(\catA,\catB,R,A,C,η)$ is a pair $(♭,\asharpis)$, in which $♭$ is an
operation ``that for each $D$ takes each $g$ to an $f$'' and
$\psharpis$ is a shorthand for this proposition:
$$\begin{array}{rl}
  (∀D.∀f.\; ♭\,D\,(♯\,D\,f) = f)\, & ∧ \\
  (∀D.∀g.\; ♯\,D\,(♭\,D\,g) = g).
  \end{array}
$$
The component $\asharpis$ of the universalness is a witness that
guarantess that this proposition holds.

The types of $♭$ and $\asharpis$ are:
$$\begin{array}{rcl}
  ♭ &:& (D : \Objs_\catB) → (\catA(A,RD) → \catB(C,D)) \\[2.5pt]
  \asharpis &:& ⟦ \;         (∀D.∀f.\; ♭\,D\,(♯\,D\,f) = f) \;\; ∧ \\
             && \phantom{m} (∀D.∀g.\; ♯\,D\,(♭\,D\,g) = g) \; ⟧
  \end{array}
$$

We can use a trick with indefinite articles to obtain the type of $♭$.
Let's overload the notations $\AllProofsOf{·}$ and $\AProofOf{·}$:
with their new meanings `$\AllProofsOf{α}$' will be pronounced ``the
type of $α$'', and `$\AProofOf{α}$' will be ``an $α$'', ``an object
with the same type of $α$'', or ``an element of $\AllProofsOf{α}$''.
Then
\begin{quote}
$♭$ is an operation that for each $D$ takes each $g$ to an $f$
\end{quote}
becomes:
$$♭ \;\;=\;\; λD.λg.\AProofOf{f}
$$

The indefinite article in this $\AProofOf{f}$ is contagious: we read
the equality above not as ``$♭$ is {\sl the} operation that takes each
$D$...'' but as ``$♭$ is {\sl an} operation that...''. We don't know
the value of $λD.λg.\AProofOf{f}$ but we can calculate its type:
$$\begin{array}{rcl}
               f  &:& C→D \\
  \AllProofsOf{f} &:& \catB(C,D) \\
  \AProofOf   {f} &:& \catB(C,D) \\
               g  &:& A→RD \\
  \AllProofsOf{g} &:& \catA(A,RD) \\
                D &:& \Objs_\catB \\
  \AllProofsOf{λD.λg.\AProofOf{f}} &:& \Objs_\catB → (\catA(A,RD) → \catB(C,D)) \\
  \end{array}
$$

We will see how to represent universalness in diagrams in
sec.\ref{ness}.

\subsection{``Physicists' notation''}

Books like \cite{ThompsonGardner} use a notation in which expressions
like ``$z=z(x,y(x))$'' are allowed --- the same symbol can be used
both as a (dependent!) variable and as the name of a function, and
arguments can be omitted --- and in a context in which $y=y(x)$ the
default meaning for $y_0$ is $y(x_0)$. In many areas of Mathematics
that notation has been phased out (see \cite[sec.3.3]{DSLsofMath}) and
replaced by one in which the names of the bound variables matter very
little.

Let's call the older notation ``physicist's notation'' and the newer
one ``mathematician's notation''; these names are not standard at all.
If we apply the ideas of the ``physicist's notation'' to judgments we
can abbreviate
$$a:A, b:B(a), c:C(a,b) ⊢ d(a,b,c):D(a,b,c)
$$
to just $a:A, b:B, c:C ⊢ d:D$, or even to $a,b,c⊢d$. Some of the
conventions of DL were inspired by conventions from ``physicist's
notation''.

\newpage

%
%
\section{The Basic Example as a skeleton \DONE}
\label{basic-example-as-skel}

In the sections \ref{the-conventions} and \ref{to-deserve-a-name} I
claimed that the diagram of the Basic Example is a ``skeleton'' of a
certain theorem, in the sense that both the statement and the proof of
that theorem can be reconstructed from just the diagram and very few
extra hints. Let's see the details of this.

%
%
\subsection{Reconstructing its functors \DONE}
\label{basic-example-functors}

Let's call this diagram --- the diagram of the Basic Example ---
$\Yzero$:
$$\Yzero \qquad := \quad \diag{Basic-Example}$$

We don't know yet the precise meaning of the functors $\catB(C,-)$ and
$\catA(A,R-)$, but if we enlarge $\Yzero$ to this --- note that we are
omitting the curved bijection by convenience,
%
%
$$\pu
  \Yzeroplus \qquad := \quad \slowdiag{Basic-Example-plus}
$$
and we draw the internal views of $\catB(C,-)$ and $\catA(A,R-)$, then
the meanings for $\catB(C,-)$ and $\catA(A,R-)$ become obvious:
%
$$\pu
  \diag{basic-example-obvious-functors}
$$

So:
$$\begin{array}{rcl}
  \catB(C,-)   &:&  \catB → \Set \\
  \catB(C,-)_0 &:=& λD.\catB(C,D) \\
  \catB(C,-)_1 &:=& λg.λf.g∘f \\
  [5pt]
  \catA(A,R-)   &:&  \catB → \Set \\
  \catA(A,R-)_0 &:=& λD.\catA(A,RD) \\
  \catA(A,R-)_1 &:=& λg.λh.Rg∘h \\
  \end{array}
$$

%

\subsection{Natural transformations}
\label{basic-example-NTs}

%
\pu

\pu

\sa{Y72+l-basicnames}{{
              \sa{A1}{A}
  \sa{A2}{C}  \sa{A3}{RC}
  \sa{A4}{D}  \sa{A5}{RD}
  \sa{A6}{E}  \sa{A7}{RE}
  \sa{B0}{\catB}      \sa{B1}{\catA}
  \sa{C0}{\catB(C,-)} \sa{C1}{\catA(A,R-)}
              \sa{A13}{η}
  \sa{A24}{f} \sa{A35}{Rf}
  \sa{A46}{g} \sa{A57}{Rg}
  \sa{A15}{\sm{ℓ,\\αDf}}
  \sa{B01}{R}
  \sa{C01}{α}
  \diag{Y72+l}
  }}

\sa{Y72+klm-idC}{{
              \sa{A1}{A}
  \sa{A2}{C}  \sa{A3}{RC}
  \sa{A4}{C}  \sa{A5}{RC}
  \sa{A6}{D}  \sa{A7}{RD}
  \sa{B0}{\catB}      \sa{B1}{\catA}
  \sa{C0}{\catB(C,-)} \sa{C1}{\catA(A,R-)}
              \sa{A13}{η}
  \sa{A24}{\sm{ι,\\\id_C}} \sa{A35}{Rι}
  \sa{A46}{f} \sa{A57}{Rf}
  \sa{A26}{\sm{f∘ι,\\f∘\id_C,\\f}}
  \sa{A15}{\sm{m,\\αCι,\\αC\id_C}}
  \sa{A17}{\sm{Rf∘(αCι),\\αD(f∘ι),\\Rf∘(αC\id_C),\\αDf}}
  \sa{B01}{R}
  \sa{C01}{α}
  \diag{Y72+klm}
  }}


%
\pu

\sa{Y0-NT-sqcond-alpha}{{
  \sa{A0}{D} \sa{B0}{\catB(C,D)} \sa{B1}{\catA(A,RD)}
  \sa{A1}{E} \sa{B2}{\catB(C,E)} \sa{B3}{\catA(A,RE)}
             \sa{C0}{\catB(C,-)} \sa{C1}{\catA(A,R-)}
  \sa{D0}{f}   \sa{D1}{αDf}         \sa{E1}{ℓ}
               \sa{D3'}{Rg∘(αDf)}   \sa{E3'}{Rg∘ℓ}
  \sa{D2}{g∘f} \sa{D3}{αE(g∘f)}
  \sa{A01}{g}
  \sa{B01}{} \sa{B02}{} \sa{B13}{} \sa{B23}{}
  \sa{C01}{α}
  \diag{Y0-NT-sqcond-inner}
}}

\sa{Y0-NT-sqcond-eta}{{
  \sa{A0}{D} \sa{B0}{\catB(C,D)} \sa{B1}{\catA(A,RD)}
  \sa{A1}{E} \sa{B2}{\catB(C,E)} \sa{B3}{\catA(A,RE)}
             \sa{C0}{\catB(C,-)} \sa{C1}{\catA(A,R-)}
  \sa{D0}{f}   \sa{D1}{Rf∘η}         \sa{E1}{ℓ}
               \sa{D3'}{Rg∘(Rf∘η)}   \sa{E3'}{Rg∘ℓ}
  \sa{D2}{g∘f} \sa{D3}{R(g∘f)∘η}
  \sa{A01}{g}
  \sa{B01}{} \sa{B02}{} \sa{B13}{} \sa{B23}{}
  \sa{C01}{α_η}
  \diag{Y0-NT-sqcond-inner}
}}

\sa{Y0-NT-sqcond-idC-iota}{{
  \sa{A0}{C} \sa{B0}{\catB(C,C)} \sa{B1}{\catA(A,RC)}
  \sa{A1}{D} \sa{B2}{\catB(C,D)} \sa{B3}{\catA(A,RD)}
             \sa{C0}{\catB(C,-)} \sa{C1}{\catA(A,R-)}
  \sa{D0}{ι}   \sa{D1}{αCι}          \sa{E1}{m} 
               \sa{D3'}{Rf∘(αCι)}    \sa{E3'}{Rf∘m}
  \sa{D2}{f∘ι} \sa{D3}{αD(f∘ι)}
  \sa{A01}{f}
  \sa{B01}{} \sa{B02}{} \sa{B13}{} \sa{B23}{}
  \sa{C01}{α}
  \diag{Y0-NT-sqcond-inner}
}}

\sa{Y0-NT-sqcond-idC-id}{{
  \sa{A0}{C} \sa{B0}{\catB(C,C)} \sa{B1}{\catA(A,RC)}
  \sa{A1}{D} \sa{B2}{\catB(C,D)} \sa{B3}{\catA(A,RD)}
             \sa{C0}{\catB(C,-)} \sa{C1}{\catA(A,R-)}
  \sa{D0}{\id_C} \sa{D1}{αC\id_C}          \sa{E1}{m} 
                 \sa{D3'}{Rf∘(αC\id_C)}    \sa{E3'}{Rf∘m}
  \sa{D2}{f}     \sa{D3}{αDf}
  \sa{A01}{f}
  \sa{B01}{} \sa{B02}{} \sa{B13}{} \sa{B23}{}
  \sa{C01}{α}
  \diag{Y0-NT-sqcond-inner}
}}


In sec.\ref{internal-view-NT} we saw that a natural transformation is
a pair. An NT $α:\catB(C,-)→\catA(A,R-)$ is a pair $(α_0, \sqcond_α)$,
where $\sqcond_α$ is this:
$$\begin{array}{rcl}
  \sqcond_α &=& ∀D.∀E.∀g.∀f.\; (\catA(A,Rf)∘αD) = (αE∘\catB(C,f)) \\
            &=& ∀D.∀E.∀g.∀f.\; (Rg∘(αDf)) = (αE(g∘f)) \\
            &=& ∀D.∀E.∀g.∀f.\; (Rg∘(α_0Df)) = (α_0E(g∘f)) \\
  \end{array}
$$

We can visualize what this ``means'' using the two diagrams at the top
in the next page.

\msk

Suppose that we define a natural transformation $α_η$ by saying that
$(α_η)_0 = λD.λf.Rf∘η$. Then we can either affirm that $\sqcond_{α_η}$
``is obvious'' or verify that it holds. Substituing $α_0$ by
$λD.λf.Rf∘η$ we obtain:
$$\begin{array}{rcl}
  \sqcond_{α_η} &=& ∀D.∀E.∀g.∀f.\; (Rg∘(Rf∘η)) = (R(g∘f)∘η) \\
  \end{array}
$$
which is clearly true, so $\sqcond_{α_η}$ holds, and $α_η$ is a
natural transformation for every $η:A→RC$. We can define an operation
$(η→α)$ by:
$$\begin{array}{rcl}
  (η↦α) &:=& λη.((α_η)_0, \sqcond_{α_η}) \\
  \end{array}
$$

Without abbreviations this definition would be very big. The lower
third of the diagram in the next page shows how visualize what
$\sqcond_{α_η}$ means.

\newpage

$$\scalebox{1.0}{$
  \begin{array}{l}
    \phantom{mmmm}
    \ga{Y72+l-basicnames}
    \\ \\ \\ \\
    \ga{Y0-NT-sqcond-alpha}
    \\ \\ \\
    \ga{Y0-NT-sqcond-eta}
  \end{array}
  $}
$$

\newpage

We can define an operation $(α↦η)$ by:
$$\begin{array}{lcl}
  (α↦η)   &:=& λα.\, α_0 \, C \, \id_C, \\
  (η↦α_0) &:=& λη.λD.λf.\,Rf∘η, \\
  (η↦α)   &:=& λη.((η↦α_0)(η), \; \sqcond_{\text{(something)}}). \\
  \end{array}
$$

We can prove that the operations $(η↦α)$ and $(α↦η)$ are mutually
inverse, but this is tricky. The proof contains a step that is hard to
visualize, and that is often stated as a slogan, like this (from
\cite[p.97]{Leinster} and \cite[p.61]{CWM2}):
%
%
\begin{quote}
  A natural transformation $α:\catB(C,-)→\catA(A,R-)$ \\
  is determined by its value at $\id_C$.
\end{quote}

The proof of that step requires instantiating $\sqcond_α$, i.e.,
$$\begin{array}{c}
  ∀D.∀E.∀g:D→E.\;∀f:C→D.\; (Rg∘(α_0Df)) = (α_0E(g∘f)) \\
  \end{array}
$$
at $D:=C$, $E:=D$, $g:=f$, and $f:=\id_C$. If we do this in two
sub-steps --- first $D:=C$ and $E:=D$, and then $g:=f$ and
$f:=\id_C$ --- we see that after the first sub-step we get this:
$$\begin{array}{c}
  ∀g:C→D.\;∀f:C→C.\; (Rg∘(α_0Cf)) = (α_0D(g∘f)) \\
  \end{array}
$$

The variables $g$ and $f$ have sort of changed their types, and some
people (like me!) would prefer to rename them, to, say:
$$\begin{array}{c}
  ∀f:C→D.\;∀ι:C→C.\; (Rf∘(α_0Cι)) = (α_0D(f∘ι)) \\
  \end{array}
$$

The diagrams in the next page show the renamed version.

To prove that our operations $(α↦η)$ and $(η↦α)$ are mutually inverse
we need to prove that the round trips $(α↦η↦α)$ and $(η↦α↦η)$ are both
identities. To prove that $(α↦η↦α)=\id$, let's define
$η_α := (α→η)(α)$ and $α_{η_α} := (η↦α)(η_α)$ and . The proof of
$(α↦η↦α)=\id$ includes this sequence of equalities:
$$\begin{array}{rcl}
  (α_{η_α})_0 D f &=& (η↦α_0) ((α↦η) (α)) D f \\
                  &=& (η↦α_0) (α C \id_C) D f \\
                  &=& (λη.λD.λf.\;Rf∘η) (α C \id_C) D f \\
                  &=& (λD.λf.\;Rf∘(α C \id_C)) D f \\
                  &=& Rf∘(α C \id_C) \\
                  &=& αDf \\
  \end{array}
$$
that uses our hard step in its last equality. The details, including
the proof of $(η↦α↦η)=\id$, can be found in \cite{MissingAgda}.

$$\scalebox{1.0}{$
  \begin{array}{l}
    \phantom{i.}
    \ga{Y72+klm-idC}
    \\ \\
    \ga{Y0-NT-sqcond-idC-iota}
    \\ \\
    \ga{Y0-NT-sqcond-idC-id}
  \end{array}
  $}
$$

\subsection{The full reconstruction \DONE}
\label{basic-example-full}

We have just reconstructed all the typings and definitions for the
diagram $\Yzero$. Here is the full reconstruction, except for the
``proof terms'' like $\respids$, $\assoc$, $\idL$ and $\idR$ for each
functor, $\sqcond$ for each natural transformations, and the proofs
that both round trips in the bijections are identity maps:
$$\slowdiag{Basic-Example}
  \qquad
  \begin{array}{rl}
    & \catA \text{ is a category}, \\
    & \catB \text{ is a category}, \\
    & R:\catB \to \catA, \\
    & A ∈ \catA, \\
    & C ∈ \catB, \\
    & η:A→RD, \\
      [5pt]
    & \catB(C,-)    :  \catB → \Set,   \\
    & \catB(C,-)_0  := λD.\catB(C,D),  \\
    & \catB(C,-)_1  := λg.λf.\;g∘f,      \\
      [5pt]
    & \catA(A,R-)   : \catA → \Set,    \\
    & \catA(A,R-)_0 := λD.\catA(A,RD), \\
    & \catA(A,R-)_1 := λg.λh.Rg∘h,     \\
      [5pt]
    & α : \catB(C,-) → \catA(A,R-), \\
      [5pt]
    & (η↦α_0) := λη.λD.λf.\;Rf∘η, \\
    & (α↦η)   := λα.\; αC(\id_C), \\
      [5pt]
    & \text{or:} \\
    & α_0 := λD.λf.Rf∘η, \\
    & η := αC(\id_C). \\
  \end{array}
$$


It is quite short --- {\sl if we treat the proof terms as
  ``obvious''.}

%
%
%

\newpage

%
%
\section{Extensions to the diagrammatic language \DONE}
\label{extensions}

Our diagrammatic language and the list of conventions in Section
\ref{the-conventions} can be extended --- ``by the user'' --- in
zillions of ways. Let's see some examples of extensions.

%
%
\subsection{A way to define new categories \DONE}
\label{comma-categories}

We saw in the sections \ref{internal-view-functor} and
\ref{basic-example-functors} how to use diagrams to define functors,
and in sections \ref{internal-view-NT} and \ref{basic-example-NTs} how
to define natural transformations. We can define new categories by
diagrams, too.


\pu
\def\commaobj#1#2#3#4{{
  \left(        \def\A{#1}
                \def\f{#4}
    \def\B{#2} \def\FB{#3}
     \diag{comma-obj-0}
  \right)
  }}


\def\dnAR{{(A{↓}R)}}

$$\pu
  \diag{defining-a-comma-category}
$$

\def\AProofOf   #1{\llangle#1\rrangle}
\def\AllProofsOf#1{\llbracket#1\rrbracket}

My favorite way --- a syntax sugar! --- of visualizing the comma
category $\dnAR$ is the middle third of the diagram above, in which
the objects of $\dnAR$ are depicted as L-shaped diagrams. To
understand the typings and the commutativity conditions we have to
look at the left third --- it indicates that $f$ must obey $Rf∘η=g$.
The right third shows a generic morphism in $\dnAR$ without the syntax
sugar, but we still have to look at the left third to type it. We
have:
$$\begin{array}{rl}
  \text{In a context in which}
    & \catA \text{ is a category}, \\
    & \catB \text{ is a category}, \\
    & R : \catB → \catA, \\
    & A \text{ is an object of $\catA$}, \\
  \text{we define the category}
    & \dnAR \text{ as follows:} \\
  [5pt]
  \text{An object of}
    & \dnAR \\
  \text{is a pair}
    & (C,η) \\
  \text{in which}
    & C : \catB_0 \\
  \text{and}
    & η : \Hom_\catA(A,RC); \\
  \text{so}
    & (C,η) : (C⠆\catB_0) × \Hom_\catA(A,RC) \\
  \text{and}
    & \dnAR_0 := (C⠆\catB_0) × \Hom_\catA(A,RC). \\
  [5pt]
  \text{A morphism}
    & f: (C,η) → (D,g) \text{ in $\dnAR$} \\
  \text{is an}
    & f: \Hom_\catB(C,D) \text{ such that $Rf∘η=g$}, \\
  \text{or equivalently a pair}
    & (f,\AProofOf{Rf∘η=g}); \\
  \text{we have}
    & (f,\AProofOf{Rf∘η=g}) : (f⠆\Hom_\catB(C,D))×\AllProofsOf{Rf∘η=g}, \\
  \text{so}
    & \Hom_\dnAR((C,η),(D,g)) := \\
    & (f⠆\Hom_\catB(C,D)) × \AllProofsOf{Rf∘η=g}.
  \end{array}
$$


\msk

This defines formally the first two components of the category
$\dnAR$. Remember that a category $\catC$ has seven components:
$$\catC = (\catC_0, \Hom_\catC, \id_\catC, ∘_\catC;
   \assoc_\catC, \idL_\catC, \idR_\catC)
$$
We are pretending that the other components of $\dnAR$ are ``obvious''
in the sense of Section \ref{to-deserve-a-name}. Note the we used the
notations for dependent types and witnesses of the sections
\ref{dependent-types} and \ref{witnesses}.


%
\subsection{Universalness as something extra \DONE}
\label{ness}

We can consider that an universal arrow is an arrow $η:A→RC$ with
something extra. We saw how to represent this ``something extra'' in
Type Theory: a universal arrow $η$ is a pair $(η,\univ_η)$, where
$\univ_η$ is its ``universalness'', that we defined in one way in
sec.\ref{freyd-with-functors} and in another way in
sec.\ref{indefinite-articles}.

Universalness is just one `-ness' among many. Several of these
``-ness''es have standard graphical representations: for example
pullbackness is indicated by a `$\pbsymbol7$', and monicness is
indicated by a tail like this: `$\monicto$'. \cite{FreydScedrov}
defines lots of graphical representations for ``-ness''es starting on
its page 37. We will sometimes use an `$:=$' to define a new
annotation that is an abbreviation for extra structure:
%
$$\pu
  \diag{universalness}
$$

This is pullbackness:
%
$$\pu
  \diag{pullbackness}
$$

%
\subsection{Opposite categories \DONE}
\label{opposite-categories}



\def\Aop{{\catA^\op}}

Suppose that we have a diagram $A \ton{f} B \ton{g} C$ in a category
$\catA$. There are several different notations for the corresponding
diagram in $\Aop$: for example, in \cite[p.33]{CWM2} it would be
written as $A \otn{f^\op} B \otn{g^\op} C$, while in
\cite[p.15]{AbramskyTzevelekos} as $A \otn{f} B \otn{g} C$. The
convention (COT) says that the notation in our diagrams should be as
close as possible to the notation in the original text --- so let's
see how to support the notation in \cite{AbramskyTzevelekos}, that
looks a bit harder than the one in \cite{CWM2}.

We want to define a new category, $\Aop$, using tricks similar to the
ones in Section \ref{comma-categories}, but now we can't pretend that
the new composition is obvious. We will define $(\Aop)_0$,
$\Hom_\Aop$, $\id_\Aop$, and $∘_\Aop$ without any textual
explanations, with just the diagrams to convince the readed that our
definitions are reasonable.
%
$$\pu
  \diag{A-and-Aop}
  \quad
  \begin{array}{c}
  \catA_0 =: (\Aop)_0 \\
  \\
  \Hom_\catA(A,B) =: \Hom_\Aop(B,A) \\
  \\
  \id_\catA(A) =: \id_\Aop(A) \\
  \\
  g ∘_\catA f =: f ∘_\Aop g \\
  \\
  \\
  \\
  \\
  \\
  \end{array}
$$

In the diagram below $F:\catA^\op→\catB$ is a contravariant functor,
and the $\catA$ above $\catA^\op$ indicates that $g:C→D$ is a morphism
of $\catA$, not of $\catA^\op$. I am not very happy with this trick
but I haven't found a better alternative yet.
%
$$\pu
  \diag{contravariant-functor}
$$

%
%
\subsection{The Yoneda Lemma \DONE}
\label{yoneda-lemma}

The formalization of $\Yzero$ as a series of typings and definitions
in Section \ref{basic-example-full} {\sl suggests} that {\sl some}
operations from Type Theory that can be applied on the formalization
side should be translatable to the diagram side; for example,
substitution. This one clearly works: if we substitute $\catA$ by
$\Set$ and $A$ by the set 1 we get this,
%
%
\pu
$$
  \Yzero \bmat{\catA := \Set \\ A := 1}
  \qquad = \quad
  \slowdiag{Basic-Example-using-Set-and-1}
$$

For each object $S$ of $\Set$ we have a bijection between elements of
$S$ and morphisms $1→S$. We will denote the morphism from 1 to $S$
that ``chooses'' an element $s∈S$ by $\nameof{s}$; the pronounciation
of $\nameof{s}$ is ``the name of $s$''. We have a bijection between
`$s$'s and $\nameof{s}$s:

For each $D∈\catB$ we have a bijection $\Set(1,RD) ↔ RD$ --- and we
can use these bijections to build a natural isomorphism
$\Set(1,R-) ↔ R$. We will draw it vertically, and complete the
triangle:

%
\pu
$$
  \Yone
  \qquad := \qquad
  \slowdiag{Basic-Example-using-Set-and-1-and-R}
$$

We can use that natural isomorphism to obtain `$β$'s from `$α$'s, and
vice-versa, by composition. We could draw an arrow for the bijection
between `$α$'s and `$β$'s and another arrow for the bijection between
`$η$'s and elements of $RC$, but we will prefer to omit them.

We will call the diagram above $\Yone$. It doesn't have a single
top-level arrow, so we can't apply the convention (CTL) to it, and we
need to specify its ``meaning'' explicitly. We will consider that its
three bijections are top-level objects, and so the diagram $\Yone$
says that we have these bijections:
$$\begin{array}{l}
  RC \\
  ↔ \; \Set(1,RC) \\
  ↔ \; \Set^\catB(\catB(C,-),\Set(1,R-)) \\
  ↔ \; \Set^\catB(\catB(C,-),R) \\
  \end{array}
$$

The Yoneda Lemma ``is'' the bijection $RC ↔ \Set^\catB(\catB(C,-),R)$
--- check how it is defined in \cite[thm.4.2.1]{Leinster},
\cite[thm.2.2.4]{Riehl}, \cite[p.61]{CWM2}, \cite[lemma 8.2]{Awodey}.
Some books show how to calculate the element of $RC$ associated to a
given $β$ and vice-versa, and most treat $\Set(1,R-)$ as something
secondary. If we represent this idea of the Yoneda Lemma in the same
format the we used in sec.\ref{basic-example-full}, we get this:
%
%
\pu
$$\Ytwo
  \;\; := \;\;
  \slowdiag{Y1-minus}
  \qquad
  \begin{array}{rl}
    & \catB \text{ is a category}, \\
    & R:\catB \to \Set, \\
    & C ∈ \catB, \\
      [5pt]
    & \catB(C,-)    :  \catB → \Set,   \\
    & \catB(C,-)_0  := λD.\catB(C,D),  \\
    & \catB(C,-)_1  := λg.λf.\;g∘f,      \\
      [5pt]
    & e ∈ RC, \\
    & \nameof{e} := λ{*}.e, \\
    & β : \catB(C,-) → R,   \\
      [5pt]
    & (e↦β_0) := λe.λD.λf.\;Rf(e), \\
    & (β↦e)   := λβ.βC\id_C, \\
  \end{array}
$$

%

%
%
\subsection{The Yoneda embedding \DONE}
\label{yoneda-embedding}

Let's define $\Ythree$ as the result of this substituion:
$$\Ythree \;\;=\;\; \Ytwo
  \bmat{
    R := \catB(B,-) \\
    e := φ \\
  }
$$

We have:
%
%
\pu
$$\scalebox{0.9}{$
  \Ythree
  \;\; = \;\;
  \slowdiag{Y3}
  \qquad
  \begin{array}{rl}
    & \catB \text{ is a category}, \\
    & \catB(B,-):\catB \to \Set, \\
    & C ∈ \catB, \\
      [5pt]
    & \catB(C,-)    :  \catB → \Set,   \\
    & \catB(C,-)_0  := λD.\catB(C,D),  \\
    & \catB(C,-)_1  := λg.λf.\;g∘f,      \\
      [5pt]
    & φ ∈ \catB(B,C), \\
    & \nameof{φ} := λ{*}.φ, \\
    & β : \catB(C,-) → \catB(B,-),   \\
      [5pt]
    & (φ↦β_0) := λφ.λD.λf.\;\catB(B,-)_1(f)(φ), \\
    & (β↦φ)   := λβ.βC\id_C, \\
  \end{array}
  $}
$$

The formulas in $(φ↦β_0)$ and $(β↦φ)$ can be simplified, and we can
derive this other diagram from it:
%
$$\pu
  \diag{Y-emb-2022mar01}
  \qquad
  \begin{array}{rl}
    & \catB \text{ is a category}, \\
    & B ∈ \catB, \\
    & C ∈ \catB, \\
      [5pt]
    & φ : B→C, \\
    & β : \catB(C,-)→\catB(B,-), \\
    & (φ↦β_0) := λφ.λD.λf.\;f∘φ, \\
    & (β↦φ)   := λβ.βC\id_C, \\
  \end{array}
$$

Let's call it $\Yfour$. This is one of the {\sl Yoneda Embeddings};
compare that diagram with the ones in \cite[p.60]{Riehl}. In
\cite{MissingAgda} we show how to formalize it in Agda as a corollary
of $\Yzero$, $\Yone$, $\Ytwo$, and $\Ythree$.
%

Note that the functor $C \mapsto \catB(C,-)$ in $\Yfour$ is
contravariant, and we used the trick from
sec.\ref{opposite-categories} to indicate this.


\subsection{Representable functors}
\label{representable-functors}

%
Some books, like \cite{Leinster} and \cite{Riehl}, present
representable functors first, then lots of examples, and only then
they present the Yoneda Lemma. In our diagram $\Yone$ a {\sl
  representation} for the functor $R$ is a pair $(C,β)$ such that the
natural transformation $β$ is a natural iso. It is easy to see that we
have these bijections, or bi-implications:
$$\begin{array}{l}
  \text{natural iso-ness of $β$} \\
  ↔ \; \text{natural iso-ness of $α$} \\
  ↔ \; \text{universalness of $η$} \\
  \end{array}
$$
The last one holds in the diagram $\Yzero$ too.

\msk

Many of the textbook examples of representable functors are
consequences of a simple theorem that shows that every functor
$R:\catB→\Set$ with a left adjoint is representable. The diagram
$\Yfive$ below is a skeleton for a proof of that theorem:
%
%
\pu
$$\Yfive \;\; := \;\;\;\;\;
  \slowdiag{Y5}
$$

The unit arrow $η_1:1→RL1$ is universal --- note the `$\univ$' in the
upper right arrow --- and so its associated `$β$' is a natural iso; we
drew it with a `$↔$'.

Statements like ``the functor {\sl blah} is representable'' are very
common in CT texts, and their natural isos is usually written just
like ``$R \cong \catB(L1,-)$'', without a name. A good way to draw the
missing diagrams for a text should not force us to invent names, and
so we should be allowed to omit the names of arrows when this is
convenient.



This is the example (iv) in \cite[p.52]{Riehl}:

\begin{quotation}
\begin{enumerate}

\item[(iv)] The functor $U:\Ring→\Set$ is represented by the unital
  ring $\Z[x]$, the polynomial ring in one variable with integer
  coefficients. A unital ring homomorphism $Z[x]→R$ is uniquely
  determined by the image of $x$; put another way, $\Z[x]$ is the {\sl
    free unital ring on a single generator}.

\end{enumerate}
\end{quotation}

The diagram below is a good way to visualize that:
%
%
\pu
$$\slowdiag{Y5-rings}
$$

Its left half is useful for when we want to remember how that
representation is generated by an adjunction of the form $F⊣U$, but it
can be omitted. In sec.\ref{the-conventions} I drew that diagram
without its left half.

\subsection{The 2-category of categories}
\label{2-category-of-cats}


\def\Dn#1{\Downarrow\scriptstyle #1}

Natural transformations are often drawn as `$⇒$'s in the middle of
``cells'' whose walls are functors. If $F,G:\catA→\catB$ are functors
and $T:F→G$ is natural transformation, then $\catA, \catB, F, G, T$
are drawn like this:
%
$$\pu
  \slowdiag{2-cat-1}
$$


There are several ways to compose functors and natural transformations
--- see \cite[section 1.7]{Riehl}, \cite[section A.5]{Hazratpour}, and
\cite{PowerPasting} for the details and the precise terminology. For
example, in
%
%
$$\pu
  \slowdiag{2-cat-2}
$$
we used ``whiskering'' and then ``vertical composition''.

We can use internal views to lower the level of abstraction of the
diagrams above. If we draw the images of an object $A∈\catA$ by the
functors and natural transformations we get:
%
$$\pu
  \slowdiag{2-cat-2-internal-view}
$$

Note that this process of taking a ``pasting diagram'' and looking at
its internal view --- in which many names become longer and some nodes
are duplicated --- is the opposite of what people usually do in the
literature: they usually go from pasting diagrams to string diagrams,
in which most names are omitted. See \cite{MarsdenCTUSD} and
\cite[section A.5]{Hazratpour}.


\newpage

%
\subsection{Kan extensions}
\label{kan-extensions}


Kan extensions are usually drawn using 2-cells (\cite[definition
  6.1.1]{Riehl}), but they can also be drawn as adjunctions
(\cite[proposition 6.1.5]{Riehl}, \cite[section X.3]{CWM2}). Let's see
how to draw them in both ways at the same time in a way that makes the
translation clear. Here is the diagram:
%
$$\pu
  \diag{kan-1}
  \qquad
  \quad
  \slowdiag{kan-2-cells-1}
$$
We will consider right Kan extensions only.

Fix $F:\catA→\catB$ and a category $\catC$. We have a functor
$(∘F):\catC^\catB→\catC^\catA$. Suppose that it has a right adjoint,
$(∘F)⊣\Ran_F$. For each natural transformation $H:\catA→\catC$ its
image by $\Ran_f$, $R:=\Ran_FH$, is a natural transformation
$R:\catB→\catC$. We have:
$$\begin{array}{rcl}
                         (∘F) &⊣& \Ran_H \\
  \Hom_{\catC^\catA}((∘F)-,-) &≅& \Hom_{\catC^\catB}(-,\Ran_F-), \\
  \Hom_{\catC^\catA}((∘F)G,H) &≅& \Hom_{\catC^\catB}(G,\Ran_FH), \\
  \Hom_{\catC^\catA}  (G∘F,H) &≅& \Hom_{\catC^\catB}(G,\Ran_FH), \\
    \Hom_{\catC^\catA}(GF,H)  &≅& \Hom_{\catC^\catB}(G,R) \\
  \end{array}
$$
and if we substitute $[G := R]$ in $\Hom_{\catC^\catB}(G,R)$ and we
transpose $\id_R$ to the left we obtain a morphism $T:RF→H$. The pair
$(R,H)$ obeys a certain universal property, that we will call
``Ran-ness'':
%
 $$∀G. \; ∀V. \; ∃!U. \; (T·UF)=V.$$

The usual way of defining right Kan extensions is by starting with the
functors $F:\catA→\catB$ and $H:\catA→\catC$ and then saying that a
pair $(R,T)$ is {\sl a} right Kan extension of $H$ along $F$ iff it
obeys Ran-ness; the functor $\Ran_F$ and the adjunction come later.
See \cite{Riehl}, section 6.1.

Note that we don't draw the `$∀V:GF→H$' in the right half of the
diagram --- it would overwrite the rest.

\newpage

%

\subsection{All concepts are Kan extensions}
\label{all-concepts}

Both \cite{CWM2} and \cite{Riehl} have sections called ``All concepts
are Kan extensions'' --- section X.7 in \cite{CWM2} and 6.5 in
\cite{Riehl}. Now that we have a favorite way of drawing right Kan
extensions we can use it to draw diagrams for 1) binary products in
$\Set$ are right Kan extensions, 2) limits are right Kan extensions
and 3) left adjoints are right Kan extensions. Here they are.

\def\bu{{•}}
\def\bubu{{••}}

\begin{enumerate}

\item Let $\bubu$ be the discrete category with two objects, $\bu$ be
  the discrete category with one object, and $!:\bubu→\bu$ be the
  unique functor from $\bubu$ to $\bu$. Then:
%
$$\pu
  \diag{kan-binary-prods}
  \qquad
  \quad
  \slowdiag{kan-2-cells-binary-prods}
$$

\item Let $\catI$ be a finite index category --- for example, $\catI =
  \psm{&&1 \\ &&↓ \\ 2&→&3}$ --- and let $\catC$ be a category with
  finite limits. A functor $D:\catI→\catC$ is a diagram of shape
  $\catI$ in $\catC$. Let's denote by $\bm1$ the discrete category
  with a single object --- the name `$\bm1$' is more standard than
  `$•$'. Then:
%
$$\pu
  \diag{kan-prods}
  \qquad
  \quad
  \slowdiag{kan-2-cells-prods}
$$

\newpage

\item Left adjoints are right Kan extensions. If
$$\catB \two/<-`->/<200>^L_R \catA$$
is an adjunction, then $(L,ε)$ is a right Kan extension of $\id_\catB$
along $R$. In a more compact notation, $L:=\Ran_R \id_\catB$ --- but
in this case we only know the action of $\Ran_R$ on the object
$\id_\catB$, and we don't know if this $\Ran_R$ can be extended to a
``real'' functor whose domain is the whole of $\catB^\catB$. The
diagram is:
%
$$\pu
  \diag{kan-adjoints}
  \qquad
  \quad
  \slowdiag{kan-2-cells-adjoints}
$$

To show that this works we have to prove that $∀V.∃!U.(ε·UR=V)$. We
will do that by ``inverting the equation $ε·UR=V$'':
$$\pu
  \slowdiag{kan-adjoints-solving-U}
$$
The solution in $U:=VL·Gη$.



\end{enumerate}

%
\subsection{A formula for Kan extensions}
\label{kan-formula}


\def\piAHSet    {\ton{π} \catA \ton{H} \Set}
\def\piAHSetmini{\xton{H∘π} \Set}
\def\pdiag    #1{\left(\diag{#1}\right)}

The sections X.3 of \cite{CWM2} and 6.2 of \cite{Riehl} discuss a
formula for calculating Kan extensions, that defines $\Ran_FH$ as the
functor whose action on objects is:
$$B \mapsto \Lim(B{↓}F \piAHSet),$$
and its action on morphisms is ``obvious'' in the sense of Section
\ref{to-deserve-a-name}. I found this formula totally impossible to
understand until I finally found a way to visualize what it ``meant''.

For each object $B∈\catB$ the functor $B{↓}F \piAHSetmini$ can be
regarded as a diagram in $\Set$ whose shape is the shape of the comma
category $B{↓}F$. If $\catA$ and $\catB$ are finite preorder categories
and $F$ is an inclusion then $B{↓}F$ can ``inherit its shape'' from
$\catA$; inclusions of preorders are ``toy examples'' ``for
children'', but they give us some intuition to start with, and they
can help us understand the formal version that can handle more general
cases.

These are the diagrams for $\Ran_F$ as a right adjoint --- note that
we use $\Set$ instead of $\catC$ to make things less abstract,
%
%
$$\pu
  \diag{kan-1}
  \qquad
  \begin{array}{l}
    RB = \\
    (\Ran_FH)B = \\
    \Lim(B{↓}F \piAHSet) \\
  \end{array}
$$
and here are some diagrams to help us understand the comma category
$B{↓}F$ --- in the compact notation its objects have names like
$(A,β)$, but in the more visual notation they are L-shaped diagrams:
%
\pu
\def\CommaObj#1#2#3#4{
             \def\A{#1}
                \def\f{#4}
  \def\B{#2} \def\C{#3}
  \left(\diag{CommaObj}\right)
}
\def\kancommaobj#1#2#3{\psm{ &#1 \\ #2&#3 \\}}
%
$$\pu
  \diag{kan-comma}
$$

\newpage

\def\kancommaobj#1#2#3{\psm{ &#1 \\ #2&#3 \\}}

Let's see an example.

\pu

If $\catA \ton{F} \catB$ is the inclusion $\pdiag{Kan_catA} →
\pdiag{Kan_catB}$,

\msk

then
$1'{↓}F = \pdiag{Kan_1'_down_F}$ and 
$3'{↓}F = \pdiag{Kan_3'_down_F}$,

\msk

and $(1'{↓}F \piAHSetmini) =  \pdiag{Kan_1'_H}$
and $(3'{↓}F \piAHSetmini) =  \pdiag{Kan_3'_H}$;

\msk

so $R(1') = \Lim(1'{↓}F \piAHSetmini) = H_2 ×_{H_6} H_5$,

and $R(3') = \Lim(3'{↓}F \piAHSetmini) = H_5$.

We can follow the same pattern to calculate $R(2')$, $R(4')$, $R(5')$,
$R(6')$.

The square of the adjunction becomes this, in this particular case:
%
\pu
\def\KanAShaped#1#2#3{
   \def\A{#1}
   \def\B{#2}
   \def\C{#3}
   \pdiag{Kan_A-shaped}
  }
\def\KanBShaped#1#2#3#4#5#6{
   \def\A{#1}
   \def\B{#2}
   \def\C{#3}
   \def\D{#4}
   \def\E{#5}
   \def\F{#6}
   \pdiag{Kan_B-shaped}
  }
%
$$\pu
  \diag{kanadjsquare-generic}
  \qquad
  \diag{kanadjsquare-example}
$$

\newpage

\subsection{Functors as objects \DONE}
\label{functors-as-objects}

One way to treat a diagram in $\Set$ like this
%
\pu
$$
  F \qquad := \qquad \diag{evil-presheaf}
$$
as a functor is to think that that diagram is an abbreviation --- it
is just the upper-right part of a diagram like this,
%
\pu
%
$$\pu
  \diag{evil-presheaf-as-functor}
$$
where we add the extra hint that the index category $\catK$ is exactly
the kite-shaped preorder category drawn above the ``$\catK$''.

The convention (CFSh) says that the image by a functor of a diagram is
a diagram with the same shape, so according to that convention we have
$F(1) = \{24,25\}$, $F(4→5) = (\{1\} \ton{1↦1} \{0,1\})$, and so
on; so the upper right part of the diagram above {\sl defines} $F$.

Note that the single `$↦$' above the $\catK \ton{F} \Set$ stands for
several `$↦$'s, one for each object and one for each morphism, and
note that $F$ is an object of $\Set^\catK$.

%
\subsection{Geometric morphisms for children \DONE}
\label{gms-for-children}

Let $\catA$ and $\catB$ be these preorder categories, and let
$f:\catA→\catB$ be the inclusion functor from $\catA$ to $\catB$:
$$
  A := \pshAargs2345
  \qquad
  B := \pshBargs123456
$$

The left half of the diagram below is the standard definition of a
geometric morphism $f$ from a topos $\calE$ to a topos $\calF$. A
geometric morphism $f:\calE→\calF$ is actually an adjunction $f^*⊣f_*$
plus the guarantee that $f^*:\calE \ot \calF$ preserves limits, which
is a condition slightly weaker than requiring that $f^*$ has a left
adjoint. When that left adjoint exists it is denoted by $f^!$, and we
say that $f^!⊣f^*⊣f_*$ is an {\sl essential geometric morphism}. The
only non-standard thing about the diagram at the left below is that is
contains an internal view of the adjunction $f^*⊣f_*$.
$$
  \diag{gm-for-adults}
  \qquad
  \def\LG{\pshAargs{G_2}{G_3}{G_4}{G_5}}
  \def\G {\pshBargs{G_1}{G_2}{G_3}{G_4}{G_5}{G_6}}
  \def\H {\pshAargs{H_2}{H_3}{H_4}{H_5}}
  \def\RH{\pshBargs{H_2{×_{H_4}}H_3}{H_2}{H_3}{H_4}{H_5}{1}}
  \diag{gm-for-children}
$$

The right half of the diagram is a particular case of the left half.
Its lower line, $\catA \ton{f} \catB$, does not exist in the left
half. The inclusion functor $f$ induces adjunctions $f^!⊣f^*⊣f_*$ as
this,
%
$$\pu
  \diag{essential-GM-small}
$$
where $f^*$ is easy to define and $f^!$ and $f_*$ not so much --- the
standard way to define $f^!$ and $f_*$ is by Kan extensions.

The big square in the upper part of the diagram is an internal view of
the adjunction $f^*⊣f_*$, with the functors $f^*G$, $G$, $H$, and
$f_*H$ being displayed as their internal views. We can choose the sets
$G_1, \ldots, G_6$ and the morphisms between them arbitrarily, so this
is an internal view of an arbitrary functor $G:\catB→\Set$; and the
same for $H$.

The arrow $f^*G \mapsot G$ can be read as a definition for the action
of $f^*$ on objects --- it just erases some parts of the diagram ---
and the arrow $H \mapsto f_*H$ can be read as a definition for the
action of $f_*$ on objects --- $f_*$ ``reconstructs'' $H_1$ and $H_6$
in a certain natural way. It is easy to reconstruct the actions of
$f^*$ and $f_*$ on morphisms from just what is shown, and to
reconstruct the two directions of the bijection.

The big diagram above can be used 1) to convince categorists who are
not seasoned toposophers that this diagrammatic language can make some
difficult categorical concepts more accessible, and 2) as a starting
point to generate diagrams ``for children'' for several parts of the
Elephant (\cite{Elephant1}), and even to prove new theorems on
toposes. For more on (1), see \cite{OchsLucatelli} and
\cite{OchsVGMS2018}; for (2), see \cite{MDE}.

\section{Related and unrelated work \DONE}
\label{related-and-unrelated}


The diagrammatic language that I described here seems to be unrelated
to the ones in \cite{CoeckePQP} and \cite{CoeckeNewStruP} --- that
describe {\sl lots} of diagrammatic languages --- and also unrelated
to \cite{MarsdenCTUSD}. We lower the level of abstraction --- see for
example Section \ref{2-category-of-cats} --- while they (usually)
raise it.


I've taken an approach that is the opposite of \cite{CaccamoWinskel}
and \cite{CaccamoPhD}. Cáccamo and Winskel define a derivation system
that can only construct functors, natural transformations, etc, that
obey the expected naturality conditions, while we allow some kinds of
sloppinesses, like constructing something that looks like a functor
and pretending that it is a functor when it may not be. When I started
working on this diagrammatic language I had a companion derivation
system for it; \cite[Section 14]{IDARCT} mentions it briefly, but it
doesn't show the introduction rules that create (proto)functors and
(proto)natural transformations and that allow being sloppy (``in the
syntactical world''). That derivation system was incomplete in all
senses --- it even had ``rules'' that I knew how to apply in
particular cases but I didn't know how to formalize.


Some of my excuses for allowing one to pretend that a functor is a
functor and leaving the verification to a second stage come from
\cite{ChengMorally}. I learned a {\sl lot} on how mathematicians use
intuition and diagrams from \cite{Kromer} --- \cite{KromerSlides} is a
great summary --- and \cite{Corfield}, and they have helped me to
identify which characteristics of my diagrammatic language are very
unusual and may be new, and that deserve to be presented in detail.


Many of the first ideas for my diagrammatic language appeared when I
was reading \cite{SeelyBeck}, \cite{SeelyLCCC}, \cite{SeelyPLC},
\cite{Jacobs}, and \cite{SeelyDiff} and trying to draw the ``missing
diagrams'' in those papers in both the original notation and in the
``archetypal case'' (\cite[Section 16]{IDARCT}).

Many of the later ideas appeared when I was trying to understand
sheaves using a certain approach ``for children'' (\cite{PH2}): I
learned how to draw diagrams showing a Grothendieck topology, its
corresponding Lawvere-Tierney topology, and its corresponding nucleus
{\sl in particular cases}, and I knew that {\sl there had to be a way}
(in the sense of \cite{ChengMorally}) to ``lift'' these diagrams of
the correspondences in particular cases to a diagram of the
correspondences in the general case... the details of this ``lifting''
were hard to formalize, but missing details started to become clear
when I required the diagrams to be translatable to a pseudocode that
could be translated to Agda. At this moment \cite{PH2} is still
incomplete, but some of the ideas in DL were motivated by its
conceptual holes.





\pu


\printbibliography

\GenericWarning{Success:}{Success!!!}  

\end{document}